\documentclass{amsart}

 \theoremstyle{definition}
 
 \theoremstyle{remark}
 
 \numberwithin{equation}{subsection}

\begin{document}

\title[On the classification problem for C$^*$-algebras]
      {On the classification problem for C$^*$-algebras}

\author{Arzikulov Farhodjon Nematjonovich}

\address{Andizhan State University, Andizhan, Uzbekistan}

\email{arzikulovfn@rambler.ru}


\date{June 5, 2013}

\begin{abstract}

In the given article it is introduced new notions of a
C$^*$-algebra of von Neumann type I and C$^*$-algebras of types
I$_n$, II, II$_1$, II$_\infty$ and III.

It is proved that any GCR-algebra is a C$^*$-algebra of von
Neumann type I, and a C$^*$-algebra is an NGCR-algebra if and only
if this C$^*$-algebra does not have a nonzero Abelian annihilator.

Also an analog of the theorem on decomposition of a von Neumann
algebra to subalgebras of types I, II and III is proved.

In the final part it is proved that every C$^*$-factor of von
Neumann type I is a C$^*$-algebra of type I$_n$ for some cardinal
number $n$, every simple C$^*$-algebra of type II$_1$ is finite,
every simple purely infinite C$^*$-algebra is of type III and
every W$^*$-factor of type II$_\infty$ has a simple
C$^*$-subalgebra of type II$_\infty$. Finally it is formulated a
classification theorem for C$^*$-factors.
\end{abstract}

\maketitle

\section*{Introduction}

\medskip

In the theory of operator algebras the classification theory of
von Neumann algebras was well developed because in any von Neumann
algebra there exists sufficient quantity of projections. A similar
classification for general C$^*$- algebras was not developed,
because there does not exist necessary quantity of projections in
these algebras. There exist definitions of C$^*$-algebras of type
I and GCR-algebras introduced by Dixmier and Kaplansky. It is
known that these definitions are equivalent \cite{3}.

In the given article, we generalize the notions of types I, II and
III in the case of C$^*$-algebras. The situation around this
problem is: there are notions of C$^*$-algebras of type I, purely
infinite C$^*$-algebras, finite C$^*$-algebras and properly
infinite C$^*$-algebras. But in many articles these notions are
considered for simple C$^*$-algebras. There is not a
correspondence between C$^*$-algebras of type I and von Neumann
algebras of type I. Indeed, on the one hand,  we can not apply the
definition of von Neumann algebra of type I to C$^*$-algebras. On
the other hand not any von Neumann algebra of type I is a
C$^*$-algebra of type I. There arises a question: can we
generalize the definition of a von Neumann algebra of type I for
C$^*$-algebras? In the given article we give an affirmative answer
for this question.

Also in this article it is proved that for any C$^*$-algebra $A$
there exist unique C$^*$-subalgebras $A_I$, $A_{II}$, $A_{III}$ of
$A$ such that $A_I$  is a C$^*$-algebra of von Neumann type I,
there does not exist a nonzero Abelian annihilator in the algebras
$A_{II}$ and $A_{III}$, the lattice $\mathcal{P}_{A_{II}}$ of
annihilators of $A_{II}$ is locally modular, the lattice
${\mathcal{P}_{A_{III}}}$ of annihilators of $A_{III}$ is purely
nonmodular. Moreover $A_I\oplus A_{II}\oplus A_{III}$ is a
C$^*$-subalgebra of $A$ and the annihilator of $A_I\oplus
A_{II}\oplus A_{III}$ is the set $\{0\}$, i.e. $Ann_A(A_I\oplus
A_{II}\oplus A_{III})=\{0\}$.

In the final part of the article a C$^*$-algebra of type I$_n$,
C$^*$-algebras of types II, II$_1$, II$_\infty$ and III are
introduced. Then we prove that any C$^*$-factor of von Neumann
type I is a C$^*$-algebra of type I$_n$ for some cardinal number
$n$, any simple C$^*$-algebra of type II$_1$ is finite, any simple
purely infinite C$^*$-algebra is of type III and any W$^*$-factor
of type II$_\infty$ has a simple C$^*$-subalgebra of type
II$_\infty$. At the end of the article a classification theorem
for C$^*$-factors is formulated.

\bigskip

\section{Annihilators of a C$^*$-algebra}

\medskip

Let $A$ be a unital C$^*$-algebra. Recall that $A_{sa}=\{a\in A:
a^*=a\}$ and $A=A_{sa}+iA_{sa}$, $A_{sa}\cap iA_{sa}=\{0\}$. Also
$Ann_r(S)=\{a\in A: sa=0 \,\text{for}\,\,\text{all}\, s\in S\}$,
$Ann_l(S)=\{a\in A: as=0 \,\,\text{for} \,\,\text{all}\,\, s\in
S\}$, where $S\subseteq A$.

\medskip

{\bf Lemma 1.} {\it Let $A$ be a unital C$^*$-algebra and $a,b\in
A$. Then

1) if $a\in A_+$, $b\in A_{sa}$ then the following conditions are
equivalent

(a) $ab+ba=0$

(b) $ab=0$

(c) $ba=0$;

2) if $a\in A_+$, $b\in A$ then $ab+ba=0$ if and only if
$ab=ba=0$.}

{\it Proof.} 1)  (a)$\Rightarrow$(b),(c): We have $ab=-ba$ and
$aba=-ba^2$, $-a^2b=aba$, that is $a^2b=ba^2$. Then $a^2$ and $b$
commute. There exists a maximal commutative C$^*$-subalgebra
$A_o$, containing $a^2$ and $b$. Since $a=\sqrt{a^2}$ we have
$a\in A_o$. Hence $ab=ba$ and $2ab=0$, i.e. $ab=ba=0$.

(b)$\Rightarrow$(a), (c): Now, suppose $ab=0$; then $ba=(ab)^*=0$
and $ba=0$. Hence $ab+ba=0$. The implication (c)$\Rightarrow$(a)
is also obvious.

2) Let $b=x+iy$, $x,y\in A_{sa}$. We have $ab+ba=ax+iay+xa+iya=0$
and $b^*a+ab^*=xa-iya+ax-iay=0$. Hence
$ab+ba+b^*a+ab^*=2(ax+xa)=0$, that is  $ax+xa=0$. Similarly
$ay+ya=0$. By 1) of lemma 1 $ax=xa=0$ $ay=ya=0$. Therefore
$ba=ab=0$.

Converse of the statement 2) is obvious. $\triangleright$

\medskip

Let $A$ be a C$^*$-algebra, $S\subseteq A$. Let
$Ann(S)=Ann_A(S)=\{a\in A: as+sa=0, \,\,for\,\, all\,\, s\in S\}$.
The set $Ann(S)$ we will call an {\it annihilator of the set} $S$.

Let $A$ be a unital C$^*$-algebra and $a$, $b$ be elements of the
set $A_{sa}$. Recall that, $A_+=\{a\in A_{sa}: \,\,\text{there}
\,\, \text{exists}\,\, b\in A \,\,\text{such}\,\,
\text{that}\,\,\, a=bb^*\}$. By lemma 1 for every set $S\subseteq
A_+$ we have $Ann(S)=Ann_r(S)\cap Ann_l(S)$.

Let $A$ be a C$^*$-algebra on a Hilbert space $H$. Then the weak
closure in $B(H)$ of $B\subseteq A$ we denote by $w(B)$. Let
$^dV=\{a\in A: xay+yax=0,\, \text{for} \, \text{any}\, x,y\in V\}$
for an arbitrary subset $V$ of $A$. We will set an analog of
decomposition on projections using annihilators for
C$^*$-algebras. First we prove the following useful lemma.

\medskip

{\bf Lemma 2.} {\it Let $A$ be a C$^*$-algebra.  Then for each
subset $S$ of $A_+$ the sets $Ann(S)$, $Ann(Ann(S))$ are
C$^*$-subalgebras and $xAx\subseteq Ann(S)$, $yAy\subseteq
Ann(Ann(S))$ for all elements $x\in Ann(S)$, $y\in Ann(Ann(S))$.
The set
$$
^d(Ann(Ann(S)))\cap ^d(Ann(S))
$$
is a Banach space. }

{\it Proof.} We will prove that $Ann(S)$ is a C$^*$-algebra. Let
$a$, $b\in Ann(S)$. Then by lemma 1 $s(ab)+(ab)s=(sa)b+a(bs)=0$
for every $s\in S$. Hence $ab\in Ann(S)$. Since $a$ and $b$ are
chosen arbitrarily we have $Ann(S)$ is an associative algebra.
Also $Ann(S)$ is a Banach algebra by separately uniformly
continuity of multiplication. Note that all conditions of the
definition of a C$^*$-algebra is valid for $Ann(S)$. Hence
$Ann(S)$ is a C$^*$-algebra.

By the previous part of the proof we have
$Ann(S)=Ann(S)_++Ann(S)_+$. It is obvious that
$Ann(Ann(S))\subseteq Ann(Ann(S)_+)$. Let $a\in Ann(Ann(S)_+)$. In
this case, if $s\in Ann(S)_{sa}$, then $s=x+iy$, $x$, $y\in
Ann(S)_+$ and $ax+xa=ay+ya=0$. Hence $as+sa=0$. Therefore $a\in
Ann(Ann(S)_{sa})$. So $Ann(Ann(S))= Ann(Ann(S)_+)$. Thus
$Ann(Ann(S))$ is a C$^*$-algebra.

It is clear that $^d(Ann(Ann(S)))$ and $^d(Ann(S))$ are linear
space. By separately uniformly continuity of multiplication they
are Banach spaces. Then $^d(Ann(Ann(S)))\cap ^d(Ann(S))$ is also a
Banach space.

By lemma 1 and associativity of multiplication we have
$$
xAx\subseteq Ann(S), yAy\subseteq Ann(Ann(S))
$$
for all elements $x\in Ann(S)$, $y\in Ann(Ann(S))$.
$\triangleright$

\medskip

{\bf Lemma 3.} {\it Let $A$ be a C$^*$-algebra on a Hilbert space
$H$, $w(A)$ be the weak closure of $A$ in $B(H)$. Then for every
$S\subseteq A_+$ the following conditions hold:

(a) There exist projections $f,e$ in $w(A)$ such that

(1) $w(Ann(Ann(S)))=ew(A)e$,  $w(Ann(S))=fw(A)f$ and
$w({^d(Ann(Ann(S)))}\cap {^d(Ann(S))})=ew(A)f+fw(A)e$,

(2) $Ann(S)=fw(A)f\cap A$, $Ann(Ann(S))=ew(A)e\cap A$ and
$^d(Ann(Ann(S)))\cap ^d(Ann(S))=[ew(A)f\oplus fw(A)e]\cap A$;

(b)
$$
Ann[Ann(Ann(S))\oplus [^d(Ann(Ann(S))\cap ^d(Ann(S))] \oplus
Ann(S)]=\{0\}.
$$
}

{\it Proof.} (a) Since $Ann(S)$ is a C$^*$-algebra (lemma 2) there
exists an increasing approximate identity $(u_\lambda)$ in
$Ann(S)$ such that $(\forall \lambda)\Vert u_\lambda \Vert \leq
1$, $(\forall \lambda \leq \mu )u_\lambda \leq u_\mu$ and $\Vert
u_\lambda \circ a-a\Vert \to 0$ for any $a\in Ann(S)$. We
calculate $\sup u_\lambda $ in $w(A)$. By the definition of
$(u_\lambda)$ $\Vert u_\lambda \circ u_\mu-u_\mu\Vert
\to_{\lambda} 0$ for any $\mu$. Then the net $(u_\lambda\circ
u_\mu)$ weakly converges to $u_\mu$ at $\lambda\to \infty$ for
each $\mu$. At the same time, since $(u_\lambda)$ weakly converges
to the element $\sup u_\lambda$ ($\sup$ is taken in $w(A)$), then
the net $(u_\lambda\circ u_\mu)$ weakly converges to $(\sup
u_\lambda)\circ u_\mu$ at $\lambda\to\infty$ for fixed $\mu$.
Hence $(\sup u_\lambda) \circ u_\mu=u_\mu$ for each $\mu$.
Therefore the net $((\sup u_\lambda) \circ u_\mu)$ weakly
converges to $\sup u_\mu$. Also the net $((\sup u_\lambda) \circ
u_\mu)$ weakly converges to $\sup u_\lambda \circ \sup u_\mu$.
Hence $\sup u_\mu=\sup u_\lambda \circ \sup u_\mu= [\sup
u_\mu]^2$. So, $\sup u_\mu$ is a projection in $w(A)$. Let
$g:=\sup u_\mu$.

By the definition of $(u_\lambda)$ the net $(s\circ u_\lambda)$
weakly converges to $s$ for any $s\in Ann(S)$, and, at the same
time $(s\circ u_\lambda)$ weakly converges to $g\circ s$. Hence
$g\circ s=s$ for all $s\in Ann(S)$. Let $f=\sup \{r(s): s\in
Ann(S)\}$ (in $w(A)$). Then $f\leq g$. Note that $Ann(S)\subseteq
U_f(w(A))$. Hence, $f\circ u_\lambda =u_\lambda$ for all
$\lambda$. Therefore $f\circ g=g$ and $f\geq g$. So $f=g$.

Now let $a$ be an arbitrary element in $U_f(w(A))$. Then there
exists a net $(a_\alpha)$ in $A$ weakly converging to $a$. Then
the net $(\{u_\lambda a_\alpha u_\mu\})$ weakly converges to
$\{u_\lambda a u_\mu\}$ for fixed $\lambda $ and $\mu $. It is
easy to see, that $(\{u_\lambda a u_\mu\})$ weakly converges to
$U_fa$ that belongs to $U_f(w(A))$. Since $a\in U_f(w(A))$ we have
$U_fa=a$. Hence, since the set $(\{u_\lambda a_\alpha u_\mu\})$ is
a net in $Ann(S)$ in relative to indices $\alpha$, $\lambda$ and
$\mu$ we have $w(Ann(S))=U_f(w(A))$.

Now, we take $Ann(Ann(S))$. By lemma 3 $Ann(Ann(S))$ is a
C$^*$-subalgebra of $A$. Hence there exists an increasing
approximate identity $(v_\lambda)$ in $Ann(Ann(S))$. Let $g=\sup
v_\lambda$ and $e=\sup \{r(s): s\in Ann(Ann(S))\}$ (in $w(A)$)
Then repeating of the above arguments gives us that $g$ is a
projection in $w(A)$ and $e=g$.

The proof of the second part of a): Note that $r(a)r(b)=0$ for all
$a\in Ann(S)$ and $b\in Ann(Ann(S))$, where $r(c)$ is the range
projection of $c\in w(A)$. Let $e=\sup \{r(a):a\in Ann(Ann(S))
\}$, $f=\sup \{r(b):b\in Ann(S) \}$. By the definitions of $e$,
$f$ we have $ef=0$.

Let $Ann_{w(A)}(S)$ be the annihilator of the set $S$ in $w(A)$.
Then there exists a projection $p$ in $w(A)$ such that
$Ann_{w(A)}(S)=pw(A)p$. At the same time, we have $Ann(S)\subseteq
Ann_{w(A)}(S)$ and $Ann(S)=Ann_{w(A)}(S)\cap A$. Then
$Ann(S)=pw(A)p\cap A$. Hence $f\leq p$, $Ann(S)=fw(A)f\cap A$ and
$Ann(Ann(S))=ew(A)e\cap A$. It can be straightforwardly prowed
that $^d(Ann(Ann(S)))\cap ^d(Ann(S))=[ew(A)f\oplus fw(A)e]\cap A$.
We have $xAy+yAx\subseteq ^d(Ann(Ann(S)))\cap ^d(Ann(S))$ for all
$x\in Ann(S)$ and $y\in Ann(Ann(S))$. Hence $Ann(S)\neq \{0\}$ and
$Ann(Ann(S))\neq \{0\}$ if $^d(Ann(Ann(S)))\cap ^d(Ann(S))\neq
\{0\}$, but not only if, because the case when $A=Ann(S)\oplus
Ann(Ann(S))$ may be valid.

(b) follows by the equality $Ann[Ann(Ann(S))\oplus Ann(S)]=\{0\}$.
$\triangleright$

\medskip

{\bf Corollary 4.} {\it Let $A$ be a C$^*$-algebra on a Hilbert
space $H$, $w(A)$ be the weak closure of $A$ in $B(H)$. If
$$
^d(Ann(Ann(S)))\cap ^d(Ann(S))=\{0\}
$$
for each $S\subseteq A_+$ then $Ann(S)$, $Ann(Ann(S))$ are two
sided ideals. In this case there exist central projections $f$,
$e$ in $w(A)$ such that $w(Ann(Ann(S)))=ew(A)e$,
$w(Ann(S))=fw(A)f$. }

{\it Proof.} By lemma 3 there exist projections  $f$, $e$ in
$w(A)$ such that $w(Ann(Ann(S)))=ew(A)e$ and $w(Ann(S))=fw(A)f$.
By the condition and separately weakly continuity of
multiplication we have
$$
w(^d(Ann(Ann(S)))\cap ^d(Ann(S)))=\{0\}.
$$
Let $p=e+f$. Then
$$
pw(A)p=ew(A)e\oplus fw(A)f,
$$
and $e$, $f$ are central projections in $pw(A)p$.

We assert that the map $\phi :A\to pAp$, defined as $\phi(a)=pap$,
for all $a\in A$, is a one-to-one correspondence between $pAp$ and
$A$. Indeed, let $a$, $b$ be elements of $A$. Suppose
$\phi(a)=\phi(b)$, i.e. $pap=pbp$. Let $x=a-b$, $C^*(x)$ be a
C$^*$-algebra, generated by $x$. It is clear that $pC^*(x)p=0$ by
separately uniformly continuity of multiplication. Let
$C^*(x)_{sa}=\{y\in C^*(x): y^*=y\}$ and $C^*(x)_+=\{y\in C^*(x):
y=zz^*,\,\text{for}\,\, \text{some}\,\,z\in C^*(x)\}$. Then
$$
C^*(x)=C^*(x)_{sa}+iC^*(x)_{sa}
$$
and
$$
C^*(x)_{sa}=C^*(x)_+-C^*(x)_+.
$$
We have $pyp=0$ for every $y\in C^*(x)_+$. Hence $py+yp=0$ for
every $y\in C^*(x)_+$. Therefore $y\in Ann(Ann(S)\oplus
Ann(Ann(S)))$. Since multiplication is separately uniformly
continuous we have $py+yp=0$ for every $y\in C^*(x)$ and
$C^*(x)\subseteq Ann(Ann(S)\oplus Ann(Ann(S)))$. From
$C^*(x)\subseteq A$ and $Ann(Ann(S)\oplus Ann(Ann(S)))=\{0\}$ it
follows that, if $x\neq 0$ then this is a contradiction. So, $x=0$
and $a=b$. Since $a$, $b$ are chosen arbitrarily we have the map
$\phi :A\to pAp$ is a one-to-one correspondence.

Now we prove that $Ann(S)$ is a closed two sided ideal of $A$. Let
$s$ be an arbitrary element of $Ann(S)$, $a$ be an arbitrary
element of $A$ and let $v$ be an arbitrary element of
$Ann(Ann(S))$. Then, since $p$, $e$, $f$ are central projections
in $pw(A)p$ and $psap\in pw(A)p$ we have
$$
(psap)v=pesapv=e(psap)fv=ef(psap)v=0,
$$
$$
v(psap)=vpesap=vfe(psap)=0.
$$
Hence
$$
p(sav+vsa)p=psavp+pvsap=psapvp+pvpsap=psapv+vpsap=0.
$$
Note that $sav+vsa\in A$. At the same time by the previous part of
the proof $pap=0$ for each $a\in A$ if and only if $a=0$. Hence
$sav+vsa=0$. Therefore, $sa \in Ann(Ann(Ann(S)))$ since $v$ is
chosen arbitrarily. But $Ann(Ann(Ann(S)))=Ann(S)$. Hence $sa \in
Ann(S)$. Hence, $Ann(S)A\subseteq Ann(S)$ since the elements $s$,
$a$ are chosen arbitrarily. Similarly $AAnn(S)\subseteq Ann(S)$
and $Ann(S)$ is a uniformly closed two sided ideal of $A$.
Similarly $Ann(Ann(S))$ is also a closed two sided ideal of $A$.
Then $w(Ann(S))$ and $w(Ann(Ann(S)))$ are closed two sided ideals
of $w(A)$ by separately weakly continuity of multiplication. Hence
$p$, $e$ and $f$ are central projections of $w(A)$.
$\triangleright$

\bigskip

\section{Lattice of annihilators of a C$^*$-algebra}

\medskip

Recall that a lattice $L$ with zero $\bf 0$, unit $\bf 1$ and an
one parameter operation ({\it orthocomplementation})
$(\,\cdot\,)^\perp:L\to L$ is called {\it an ortholattice} if $L$
satisfies the following conditions

\smallskip
 {\it(1)}~$ x\wedge x^\perp={\bf 0},\quad
 x\vee x^\perp={\bf 1};$

 \smallskip
 {\it(2)}~$x^{\perp\perp}\!:=(x^\perp)^\perp=x;$

 \smallskip
  {\it(3)}~$(x\vee y)^\perp=x^\perp\wedge y^\perp,\quad
 (x\wedge y)^\perp=x^\perp\vee y^\perp.$

\smallskip

An ortholattice $L$ is called {\it an orthomodular lattice\/}, if
the {\it orthomodular law\/} is valid in this lattice: for all $x,
y\in L$, from $x\le y$ follows $y=x\vee (y\wedge x^\perp)$.

Let $x$ and $y$ be elements of an ortholattice $L$. If $x=(x\wedge
y)\vee(x\wedge y^\perp)$ then we say $x$ commutes with $y$ and
write $x {\bf C} y$. It is clear that $x {\bf C} y$ if $x\leq y$.
The relation ${\bf C}$ is not a symmetric relation.

Recall that a lattice is said to be {\it modular,\/} if it follows
from $x$, $z\in L$, $x\leq z$ that $x\vee(y\wedge z)=(x\vee
y)\wedge z$ for every $y\in L$.

A subset $B$ of an orthomodular lattice $L$ is called {\it a
boolean subalgebra\/}, if $B$ is a boolean algebra with the
induced lattice operations and the orthocomplement in the sense of
boolean complement. Maximal elements of the set of all boolean
subalgebras of $L$ ordered by inclusion we call {\it maximal
boolean subalgebras\/} of $L$. By the Kuratovskiy-Zorn's lemma for
every boolean subalgebra there exists a maximal boolean subalgebra
containing this boolean subalgebra. But the following improved
result holds.

An orthomodular lattice is a boolean algebra if and only if any
two elements of this lattice are compatible.

Intersection of all maximal boolean subalgebras of an orthomodular
lattice $L$ is called a center ${\bf Z}(L)$ of the orthomodular
lattice $L$. It is clear that the center ${\bf Z}(L)$ consists of
elements compatible with all elements of $L$. The center of an
orthomodular lattice is a boolean subalgebra.

A lattice $L$ is said to be {\it order complete}, if for every
subset $M\subset L$ there exists a least upper bound $\bigvee
M\!:=\sup(M)$ in $L$. Of course, in this case, if $L$ is
orthomodular, then also there exists $\bigwedge M\!:=\inf(M)$ and
 $\bigwedge M=\Big(\bigvee_{x\in M}x^\perp\Big)^\perp$.

The center ${\bf Z}(L)$ of a complete orthomodular lattice $L$ is
a complete boolean algebra.

Let $A$ be a C$^*$-algebra. We introduce the set $\mathcal{P}$ of
all annihilators for $A$ as follows
$$
\mathcal{P}=\{V\subset A: \,\,there\,\, exists\,\, S\subseteq A_+
\,\, such\,\,that\,\,V=Ann(Ann(S))\}.
$$
Note that, since $Ann(Ann(Ann(S)))=Ann(S)$ we have
$$
\mathcal{P}=\{V\subset A: \,\,there\,\, exists\,\, S\subseteq A_+
\,\, such\,\,that \,\,V=Ann(S)\}.
$$
For every two elements $V$, $W$ of $\mathcal{P}$, if $V\subseteq
W$, then we write $V\le W$. So we define an {\it order} in
$\mathcal{P}$.

\medskip

{\bf Lemma 5.} {\it Let $A$ be a C$^*$-algebra, $\mathcal{P}$ be
the set of annihilators, defined above. Then $(\mathcal{P},\le )$
is a complete lattice. }

{\it Proof.} Let $V$, $W$ be elements in $\mathcal{P}$. Then there
exist $S$, $P\subseteq A$ such that $V=Ann(Ann(S))$,
$W=Ann(Ann(P))$ respectively. It is clear that $V$, $W\subseteq
Ann(Ann(P\cup S))$. Let $Z\in\mathcal{P}$ such that $V\subseteq
Z$, $W\subseteq Z$. Then there exists
 $Q\subseteq A_+$ such that $Ann(Ann(Q))=Z$. We note that
$Ann(Ann(Ann(Q)))\subseteq Ann(Ann(Ann(S)))$. At the same time,
$Ann(Ann(Ann(S)))=Ann(S)$ and $Ann(Ann(Ann(Q))=Ann(Q)$. Hence,
$Ann(Q)\subseteq Ann(S)$. Similarly $Ann(Q)\subseteq Ann(P)$.
Hence by the definition of an annihilator $Ann(Q)\subseteq
Ann(P\cup S)$. Therefore $Ann(Ann(P\cup S))\subseteq
Ann(Ann((Q))$. Since $Z$ is chosen arbitrarily we have $V\vee
W=Ann(Ann(P\cup S))$.

Note that $Ann(Ann(P)\cup Ann(S))\subseteq V\cap W$. Let
$Z\in\mathcal{P}$ such that $Z\subseteq V$, $Z\subseteq W$. Then
there exists $Q\subseteq A_+$ such that $Ann(Ann(Q))=Z$. By the
definition of an annihilator we have $Ann(Ann(Q))\subseteq
Ann(Ann(S)\cup Ann(P))$. Since $Z$ is chosen arbitrarily $V\wedge
W=Ann(Ann(P)\cup Ann(S))$.

Note, if $S\subseteq A_+$ then $Ann(Ann(S)\cup S)=\{0\}$,
$$
\sup \{Ann(S), Ann(Ann(S))\}=\sup \{Ann(Ann(Ann(S))),
Ann(Ann(S))\}=
$$
$$
Ann(Ann(Ann(S)\cup S))=A
$$
and $Ann(S)\wedge Ann(Ann(S))\subseteq Ann(S)\cap Ann(Ann(S))
=\{0\}$, i.e. $Ann(S)\wedge Ann(Ann(S))=\{0\}$.

Hence, the set $\mathcal{P}$, equipped with the order $\subseteq$,
is a lattice.

Let $\{V_i\}$ be an arbitrary subset of $\mathcal{P}$. Then there
exist $\{S_i\}\subseteq A_+$ such that $Ann(Ann(S_i))=V_i$  for
all $i$. We have $V_i\subseteq Ann(Ann(\cup_i S_i))$ for each $i$.
Let $Z$ be an element in $\mathcal{P}$ such that $V_i\subseteq Z$,
for each $i$. Then there exists $Q\subseteq A_+$ satisfying the
condition $Ann(Ann(Q))=Z$. Note that $Ann(Ann(Ann(Q)))\subseteq
Ann(Ann(Ann(S_i)))$ for every $i$. At the same time,
$Ann(Ann(Ann(S_i)))=Ann(S_i)$ and $Ann(Ann(Ann(Q))=Ann(Q)$. Hence,
$Ann(Q)\subseteq Ann(S_i)$ for every $i$. Hence by the definition
of an annihilator $Ann(Q)\subseteq Ann(\cup_i S_i)$. Therefore
$Ann(Ann(\cup_i S_i))\subseteq Ann(Ann((Q))$. Since $Z$ is chosen
arbitrarily $\bigvee_i V_i=Ann(Ann(\cup_i S_i))$. Hence the
lattice $(\mathcal{P},\le )$ is complete. $\triangleright$

\medskip

{\bf Lemma 6.} {\it Let $A$ be a C$^*$-algebra and $X$, $Y\in
\mathcal{P}$. Then

(a) $X\wedge Ann(X)=\{0\}$, $X \vee Ann(X)=A$;

(b) $Ann(Ann(X))=X$, and if $X\neq A$ then $Ann(X)\neq \{0\}$;

(c) $Ann(X\vee Y)=Ann(X) \wedge Ann(Y)$,
    $Ann(X \wedge Y)=Ann(X) \vee Ann(Y)$.
}

{\it Proof.} (a) Let $S$ is a subset of $A_+$ and $X=Ann(Ann(S))$.
Then by the proof of lemma 5 $Ann(S)\wedge Ann(Ann(S))=\{0\}$. We
have $Ann(X)=Ann(Ann(Ann(S)))=Ann(S)$. Then $Ann(X)\wedge
X=\{0\}$. Similarly $Ann(S)\vee Ann(Ann(S))=A$ and $Ann(X)\vee
X=A$.

(b) Suppose $X\neq A$ and $Ann(X)=\{0\}$; then $Ann(Ann(X))=A$.
But by the definition $Ann(Ann(X))=X$. This is a contradiction.
Hence $Ann(X)\neq \{0\}$.

(c) Let $Q\subseteq A_+$ and $Y=Ann(Ann(Q))$. By the proof of
lemma 5 $Ann(X\vee Y)=Ann(Ann(Ann(S\cup Q)))=Ann(S\cup Q)$. At the
same time $Ann(X)\wedge Ann(Y)=Ann(Ann(Ann(S)))\wedge
Ann(Ann(Ann(Q)))= Ann(S)\wedge Ann(Q)$. We have $Z\subseteq
Ann(S)\cap Ann(Q)$ for any $Z\in \mathcal{P}$ such that
$Z\subseteq Ann(S)$ and $Z\subseteq Ann(Q)$. At the same time
$Ann(S)\cap Ann(Q)=Ann(S\cup Q)$. Hence $Ann(S)\wedge
Ann(Q)=Ann(S\cup Q)$. Thus $Ann(X)\wedge Ann(Y)=Ann(S\cup Q)$ and
$Ann(X\vee Y)=Ann(X)\wedge Ann(Y)$.

Similarly we have $Ann(X\wedge Y)=Ann(X)\vee Ann(Y)$.
$\triangleright$

\medskip

{\it Example.} Let $\mathcal{X}$ be a compact, $\tau_\mathcal{X}$
be the topology of $\mathcal{X}$. Let $\leq$ be an order in
$\tau_\mathcal{X}$, defined as follows: if $V$, $W\in
\tau_\mathcal{X}$ and $V\subseteq W$ then $V\leq W$.

The ordered set $(\tau_\mathcal{X},\leq)$ is a lattice. Indeed,
${\bf 1}=\mathcal{X}$, ${\bf 0}=\{\oslash\}$, $V\vee W=V\cup W$,
$V\wedge W=V\cap W$ for all $V$, $W\in \tau_\mathcal{X}$.

The ordered set $(\tau_\mathcal{X},\leq)$ is a complete lattice.
Indeed, let $\{V_i\}\subseteq \tau_\mathcal{X}$. Then $\bigvee_i
V_i=\cup_i V_i$ and $\bigwedge_i V_i=\cup \{U\in \tau_\mathcal{X}:
\,\,for\,\, any \,\, i \,\, U\subseteq V_i\}$.

Moreover, $(\tau_\mathcal{X},\leq)$ is a complete boolean algebra.
Indeed, for arbitrary $V$, $W\in \tau_\mathcal{X}$ we have
$V=V_1\vee Z$, $W=W_1\vee Z$, where $V_1=V\setminus (V\cap W)$,
$W_1=W\setminus (V\cap W)$, $Z=V\cap W$ and $V_1$, $W_1$, $Z\in
\tau_\mathcal{X}$.

Let $C^c(\mathcal{X})$ be the complex commutative algebra of
continuous functions on $\mathcal{X}$. Then the lattice
$\mathcal{P}_{C^c(\mathcal{X})}$ of annihilators of
$C^c(\mathcal{X})$ is a complete boolean algebra. Moreover,
$\mathcal{P}_{C^c(\mathcal{X})}$ is order isomorphic to the
complete boolean algebra $(\tau_\mathcal{X},\leq)$, where the
isomorphism is defined by the map
$$
\Phi (X)=\{x\in \mathcal{X}: f(x)\neq 0\,\,for \,\,some
\,\,function\,\,f\in X\}, X\in \mathcal{P}_{C^c(\mathcal{X})}.
$$
Indeed, $U_f=\{x\in \mathcal{X}: f(x)\neq 0\}$, where $f\in X$, is
open in $\mathcal{X}$. Hence $U_X=\cup_{f\in X}U_f$ is also open
in $\mathcal{X}$. Since $\Phi (X)=U_X$ $\Phi (X)$ is an open set
in $\mathcal{X}$. The set $C(\Phi (X))$ of all functions $f\in
C^c(\mathcal{X})$ such that $\{x\in \mathcal{X}: f(x)\neq
0\}\subseteq \Phi (X)$ forms a commutative subalgebra of
$C^c(\mathcal{X})$. Moreover $C(\Phi (X))\in
\mathcal{P}_{C^c(\mathcal{X})}$ and $Ann(C(\Phi (X)))=Ann(X)$.
Hence $C(\Phi (X))=X$. Let $Y\in \mathcal{P}_{C^c(\mathcal{X})}$
and $\Phi (X)=\Phi (Y)$. Then $C(\Phi (Y))=Y$ and $X=Y$ by $\Phi
(X)=\Phi (Y)$.

Let $A$ be a C$^*$-algebra. An annihilator $V\in \mathcal{P}$ is
said to be {\it central}, if
$$
^d(Ann(Ann(S)))\cap {}^d(Ann(S))=0,
$$
where $S\subseteq A_+$ and $V=Ann(Ann(S))$. The set of all central
annihilators we denote by $Z(\mathcal{P})$. We will say that two
annihilators $V$ and $W$ in $\mathcal{P}$ are {\it orthogonal}, if
$V\cdot W=\{0\}$, where $V\cdot W=\{vw:v\in V, w\in W\}$.

\medskip

{\bf Lemma 7.} {\it Let $A$ be a C$^*$-algebra on a Hilbert space
$H$, $Z(\mathcal{P})$ be the set of all central annihilators in
$\mathcal{P}$. Then elements of $Z(\mathcal{P})$ are pairwise
commute, i.e. $X=(X \wedge Y)\vee (X \wedge Y^\perp)$ for any $X$,
$Y\in Z(\mathcal{P})$. }

{\it Proof.} Let $X$, $Y\in Z(\mathcal{P})$. Then $X \wedge Y=X
\cap Y$, $X \wedge Y^\perp=X \cap Y^\perp$ and by lemma 3 $X=A\cap
ew(A)e$, $Y=A\cap fw(A)f$, $Y^\perp=A\cap (\bar{f})w(A)(\bar{f})$
for some projections $e$, $f$ and $\bar{f}$ in $w(A)$. Note that
$Ann(X)\subseteq Ann((X \cap Y)\cup (X \cap Y^\perp))$ and $e$,
$f$ and $\bar{f}$ are central projections in $w(A)$. Let $X\cdot
Y=\{xy: x\in X, y\in Y\}$. Then $X\cdot Y\subset X\cap Y$, $X\cdot
Y\subset efw(A)ef$, $X\cap Y=A\cap efw(A)ef$. Similarly $X\cap
Y^\perp=A\cap e\bar{f}w(A)e\bar{f}$.

Suppose $Ann(X)\neq Ann((X \cap Y)\cup (X \cap Y^\perp))$; then
there exists $a\in Ann((X \cap Y)\cup (X \cap Y^\perp))$ such that
$a\notin Ann(X)$. Hence there exists $x\in X_+$ such that $ax\neq
0$. Since $Ann(Y\cup Y^\perp)=\{0\}$ then there exists $y\in Y\cup
Y^\perp$ such that $(ax)y\neq 0$. We have
$(ax)y(f+\bar{f})e=(ax)y$. Then $(ax)y\in
(f+\bar{f})ew(A)(f+\bar{f})e$. But $a(fe+\bar{f}e)=0$. Hence $a\in
Ann(X)$. Therefore $Ann(X)= Ann((X \cap Y)\cup (X \cap Y^\perp))$
and since $(X \cap Y)\cup (X \cap Y^\perp)\in Z(\mathcal{P})$ we
have $X=(X \wedge Y)\vee (X \wedge Y^\perp)$. $\triangleright$

\medskip

{\bf Lemma 8.} {\it Let $A$ be a C$^*$-algebra on a Hilbert space
$H$, $Z(\mathcal{P})$ be the set of all central annihilators of
$\mathcal{P}$. Then $Z(\mathcal{P})$ is a complete boolean
algebra. }

{\it Proof.} By lemma 7 elements of $Z(\mathcal{P})$ pairwise
commute. Hence $Z(\mathcal{P})$ is a boolean algebra by the first
part of this section.

Let $\{V_i\}$ be a subset of $Z(\mathcal{P})$. Then by corollary 4
for any $i$ there exist central projections $e_i$, $f_i\in
P(w(A))$ such that $w(V_i)=e_i(w(A))$, $w(Ann(V_i))=f_i(w(A))$,
where $w(S)$ is the weak closure of a set $S\subseteq A$ in
$B(H)$. Then $V_i$, $Ann(V_i)$ are uniformly closed two sided
ideals of $A$ for all indices $i$.

Let $a$, $v$ be arbitrary elements of $A$, $\cap_i Ann(V_i)$,
respectively. Then $v\in Ann(V_i)$ and $av$, $va$ belong to
$Ann(V_i)$ for all $i$. Hence $av$, $va$ belong to $\cap_i
Ann(V_i)$ to. Therefore, $\cap_i Ann(V_i)$ is a two sided
uniformly closed ideal of $A$ by $\cap_i Ann(V_i)=Ann(\cup_i
V_i)$. There exists a projection $f\in w(A)$ such that $w(\cap_i
Ann(V_i))=fw(A)f$. Then by separately weakly continuity of
multiplication the following equality is valid
$$
^d(Ann(Ann(V_i)))\bigcap {}^d(Ann(V_i))={0}.
$$
Therefore, by corollary 4 $Ann(Ann(\cup_i V_i))\in
Z(\mathcal{P})$. At the same time, by the proof of lemma 5 $\sup_i
V_i=Ann(Ann(\cup_i V_i))$. Hence $\sup_i V_i\in Z(\mathcal{P})$.

Similarly $\inf_i V_i\in Z(\mathcal{P})$. So the lattice
$Z(\mathcal{P})$ is complete. $\triangleright$

Let $V\in \mathcal{P}$. By lemma 5 the greatest lower bound $c(V)$
of central annihilators $W\in Z(\mathcal{P})$ satisfying
$V\subseteq W$, is also an annihilator. Moreover by lemma 8 $c(V)$
is central. The annihilator $c(V)$ we will call a {\it central
support} of $V$.

\medskip

{\bf Lemma 9.} {\it Let $A$ be a commutative C$^*$-algebra on a
Hilbert space $H$, $w(A)$ be the weak closure of $A$ in $B(H)$.
Let $X$ be the topological space of multiplicative functionals of
$A$, $Y$ be the topological space of multiplicative functionals of
$w(A)$. Let $supp(Y)$, $supp(X)$ be the sets of all points of the
spaces $Y$ and $X$ respectively. Then

    (a) $supp(X)\subseteq supp(Y)$,

    (b) the set $supp(X)$ of all points of the space $X$ is dense in the
topological space $Y$. }

{\bf Proof.} (a) Since every multiplicative functional on $A$ can
be uniquely $*$-weakly extended to a multiplicative functional on
the algebra $w(A)$ we may assume $supp(X)\subseteq supp(Y)$.

(b) Suppose $supp(X)$ is not dense in $Y$. Let $C(X)$, $C(Y)$ be
the commutative algebras of complex-valued continuous functions on
$X$, $Y$ respectively. Then $A\cong C(X)$, $w(A)\cong C(Y)$.

Note, that $a(x)=\bar{a}(x)$ for all $a\in C(X)$ and $x\in X$,
where $\bar{a}$ is the image of the function $a$ in $C(Y)$ in
point of $C(X)\subseteq C(Y)$. Let $Y_o$ be an open subset of $Y$
such that $Y_o\cap X=\oslash$. The set $C(Y_o)$ of all functions
$f\in C(Y)$ such that $\{x\in Y: f(x)\neq 0\}\subseteq Y_o$ forms
a commutative subalgebra of $C(Y)$ and $C(X)\subseteq
Ann_{C(Y)}(C(Y_o))$. Let $f$ be an arbitrary nonzero element of
$C(Y_o)$. Then $f\cdot C(X)=\{0\}$. By separately weakly
continuity of multiplication $f\cdot w(C(X))=f\cdot C(Y)=\{0\}$.
Hence $f=0$. This is a contradiction. Therefore $Y_o=\oslash$.
$\triangleright$

\medskip

{\bf Lemma 10.}  {\it Let $A$ be a commutative C$^*$-algebra on a
Hilbert space $H$, $\mathcal{P}$ be the set of annihilators and
$Y\in \mathcal{P}$, $X\in \mathcal{P}$. Suppose $X$ is a subset of
$Y$ such that $X\neq Y$; then $Ann_Y(X)\neq \{0\}$ and
$Ann_Y(Ann_Y(X))=X$. }

{\it Proof.} Let $w(A)$ be the weak closure of $A$ in
$B(H)$, $w(Y)$ be the weak closure of $Y$ in
$w(A)$. Then the weak closure $w(X)$ of $X$ in $w(Y)$ coincides
with $ew(Y)e$ for some projection $e\in w(Y)$ satisfying the condition
$e<1$, i.e. $w(X)=ew(Y)e$.

Let $Q$ be the topological space of multiplicative functionals of
$A$, $\bar{Q}$ be the topological space of multiplicative
functionals of $w(A)$. By (a) of lemma 9 $supp(Q)\subseteq
supp(\bar{Q})$. By (b) of lemma 9 the set $supp(Q)$ is dense in
$\bar{Q}$.

Note that $V=\{x\in \bar{Q}: e(x)\neq 0\}$ and $W=\{x\in \bar{Q}:
(1-e)(x)\neq 0\}$ are close-open subsets of $\bar{Q}$ and
$\bar{Q}=V\cup W$. Also $Q_Y=\bigcup_{f\in Y}\{x\in Q: f(x)\neq
0\}$, $Q_X=\bigcup_{f\in X}\{x\in Q: f(x)\neq 0\}$ are open
subsets of $Q$. Let $Cl(Q_Y)$ be the closure of $Q_Y$ and
$Cl(Q_X)$ be the closure of $Q_X$ in $Q$. If $Q_Y\neq Q_X$ then
$Cl(Q_Y)\neq Cl(Q_X)$. Indeed, if $Cl(Q_Y)=Cl(Q_X)$ then
$Q\setminus Cl(Q_Y)=Q\setminus Cl(Q_X)$, $Q\setminus Cl(Q_Y)$ is a
nonempty open set in $Q$ and $Ann_A(Y)=Ann_A(X)$ (see the example
above). Then $Y=X$. This is imposable. Hence $Cl(Q_Y)\neq Cl(Q_X)$
and $Q_Y\neq Q_X$. Otherwise also we get $Cl(Q_Y)=Cl(Q_X)$.

Then $Q_Y\setminus Q_X$ is an open set in $Q$ since $Q_X\subset
Q_Y$. Therefore $Ann_Y(X)\neq \{0\}$.

Since $Q_Y=(Q_Y\setminus Q_X)\cup Q_X$ we have
$Ann_Y(Ann_Y(X))=X$. The proof is completed. $\triangleright$

\medskip

{\it Designation.} Let $A$ be a C$^*$-algebra on a Hilbert space
$H$, $\mathcal{P}$ be the set of annihilators. Then by (A) we
denote the following condition:

(A) : For every annihilator $V\in \mathcal{P}$ and for every
maximal commutative $*$-subalgebra $V_o$ of $V$ the identity
element $e_{V_o}$ of $w(V_o)$ coincides with the identity element
$e_V$ of $w(V)$, i.e. $e_{V_o}=e_{V}$.

For example, each von Neumann algebra satisfies condition (A).

{\bf Lemma 11.} {\it Let $A$ be a C$^*$-algebra on a Hilbert space
$H$, $\mathcal{P}$ be the set of annihilators. Suppose $A$
satisfies condition (A). Let $Y\in \mathcal{P}$, $X\in
\mathcal{P}$. If $X$ is a subset of $Y$ such that $X\neq Y$ then
$Ann_Y(X)\neq \{0\}$ and $Ann_Y(Ann_Y(X))=X$. }

{\it Proof.} Let $w(A)$ be the weak closure of $A$ in $B(H)$,
$w(Y)$ be the weak closure of $Y$ in $w(A)$. Then by (2) of lemma
3 the weak closure $w(X)$ of $X$ in $w(Y)$ coincides with $fw(Y)f$
for some projection $f\in w(Y)$ such that $f<e$, i.e.
$w(X)=fw(Y)f$, where $e$ is an identity element of $w(Y)$.

Note that $X$ is a C$^*$-algebra. By the supposition for every
maximal commutative $*$-subalgebra $X_o$ of $X$ we have $f\in
w(X_o)$ and $w(X_o)=fw(Y_o)f$, where $Y_o$ is a maximal
commutative $*$-subalgebra of $Y$, containing $X_o$. Since $f<e$
we have $e\notin X_o$ and $X_o\neq Y_o$. Hence by lemma 10
$Ann_{Y_o}(X_o)\neq \{0\}$ and $Ann_{Y_o}(Ann_{Y_o}(X_o))=X_o$.
Since $w(X_o)=fw(Y_o)f$ and $Ann_{Y_o}(X_o)\neq \{0\}$ we have
$(e-f)Y_o\cap Y_o\neq \oslash$. Therefore $Ann_Y(X)\neq \{0\}$.

We have $Ann_{Y_o}(Ann_{Y_o}(X_o))=X_o$ for every maximal
commutative $*$-subalgebra $X_o$ of $X$ and for every maximal
commutative $*$-subalgebra $Y_o$ of $Y$, containing $X_o$. Also
$Ann_{Y_o}(X_o)\subset Ann_Y(X_o)$ and, since $f\in w(X_o)$, $f\in
w(X)$ we have $Ann_{Y_o}(X_o)\subset Ann_Y(X)$. Hence
$Ann_Y(Ann_Y(X))\subset Ann_Y(\cup_{X_o\in
Max(X)}Ann_{Y_o}(X_o))$, where $Max(X)$ is the set of all maximal
commutative $*$-subalgebras of $X$. At the same time $X_o\subset
Ann_Y(\cup_{X_o\in Max(X)}Ann_{Y_o}(X_o))$ for all $X_o\in
Max(X)$. Then  $Ann_Y(Ann_Y(X))=X$. $\triangleright$

\medskip

{\bf Theorem 12.} {\it Let $A$ be a C$^*$-algebra, $\mathcal{P}$
be the set of annihilators, defined above. Then $\mathcal{P}$ is
an ortholattice. Moreover, if the C$^*$-algebra $A$ satisfies
condition (A), then $(\mathcal{P},\le )$ is an orthomodular
lattice.}

{\it Proof.} Let $\mathcal{P}=(\mathcal{P},\le )$ and
$(\,\cdot\,)^\perp:\mathcal{P}\to \mathcal{P}$ be the map defined
as $(X)^\perp=Ann(X)$, for any $X\in \mathcal{P}$. Then by lemma 6
the map $(\,\cdot\,)^\perp$ is an orthocomplementation and
$\mathcal{P}$ is an ortholattice with this operation.

Now we prove that, if $A$ satisfies condition (A), then
$\mathcal{P}$ is orthomodular. Let $X$, $Y\in \mathcal{P}$, $S$,
$Q\subseteq A_+$, $X=Ann(Ann(S))$, $Y=Ann(Ann(Q))$ and $X\leq Y$.
Then $\inf \{Ann(Ann(Q)),Ann(S)\}=Ann_Y(Ann(Ann(S)))$. Indeed,
$$
Ann_Y(Ann(Ann(S)))=Ann(Ann(Ann(S)))\cap Ann(Ann(Q))
$$
$$
=Ann(S)\cap Ann(Ann(Q)).
$$
Hence by the proof of (c) of lemma 6
$$
\inf\{Ann(Ann(Q)),Ann(S)\}=Ann(Ann(Q))\cap Ann(S).
$$

So, $Y\wedge X^\perp=Ann_Y(Ann(Ann(S)))$ and
$Ann_Y(Ann(Ann(S)))\in \mathcal{P}$ since
$Ann_Y(Ann(Ann(S)))=Ann(Ann(Ann(S)))\cap Ann(Ann(Q))$.

Suppose $A\subseteq B(H)$ and the identity element of $B(H)$ is
the identity element of $w(A)$ for some Hilbert space $H$, where
$w(A)$ is the weak closure of $A$ in $B(H)$. Let $w(Y)$ be the
weak closure of $Y$ in $w(A)$ and $e$, $f$ be projections in
$w(A)$ such that $X=A\cap fw(A)f$, $Y=A\cap ew(A)e$. Then $f\leq
e$. If $e=f$ then the assertion of the theorem is true.

Suppose $f\neq e$. By lemma 11 $X\in \mathcal{P}_Y$. Therefore
$Ann_Y(X)=(e-f)(w(Y))(e-f)\cap Y$ and $w(Ann_Y(X))=gw(Y)g$ for
some projection $g\in w(A)$. Then $X\vee Ann_Y(X)=Y$. Hence
$Y=X\vee (Y\wedge X^\perp)$. Therefore $\mathcal{P}$ is an
orthomodular lattice. $\triangleright$

Let $A$ be a C$^*$-algebra. An annihilator $V\in \mathcal{P}$ is
said to be {\it Abelian}, if $V$ is a commutative C$^*$-subalgebra
of $A$. Let $B$ be a C$^*$-subalgebra of $A$ and
$$
\mathcal{P}_B=\{V\subseteq B: \,\, there \,\, exists \,\, such
\,\, S\subseteq B_+ \,\, that \,\, V=Ann_B(Ann_B(S))\}.
$$

\medskip

{\bf Lemma 13.} {\it Let $A$ be a C$^*$-algebra. Then the
following statements are valid.

a) Let $V\in \mathcal{P}$. Then $\mathcal{P}_V\subseteq \{W\in
\mathcal{P}: W\subseteq V\}$ and $\mathcal{P}_V$ is a complete
sublattice of $\mathcal{P}$. Moreover, if the C$^*$-algebra $A$
satisfies condition (A), then $\mathcal{P}_V=\{W\in \mathcal{P}:
W\subseteq V\}$

b) Let $V\in \mathcal{P}$ and $Z$ be a central annihilator in
$\mathcal{P}$ such that $V\subseteq Z$. Then $Ann_Z(Ann_Z(V))=V$,
i.e. $V\in \mathcal{P}_Z$, and $Ann_Z(V)=\{vzv : v\in Ann(V), z\in
Z\}$. Conversely, if for an arbitrary subset $V\subseteq Z$
$Ann_Z(Ann_Z(V))=V$, i.e. $V\in \mathcal{P}_Z$, then $V\in
\mathcal{P}$.

c) Let $V$ be an Abelian annihilator. Then for every $W\in
\mathcal{P}$, if $W\subseteq V$, then $W$ is an Abelian
annihilator to. }

{\it Proof.} a) Let $Z\in \mathcal{P}_V$. Then $Ann(V) \subseteq
Ann(Z)$ and
$$
Ann(Ann(Z))\subset Ann(Ann(V))=V.
$$
Since $Ann_V(Z)\subseteq Ann(Z)$ and
$$
Ann(Ann(Z))=Ann_V(Ann(Z))=Ann_V(Ann_V(Z))
$$
we have $Ann(Ann(Z))=Z$. Hence $Z\in \mathcal{P}$.

Now, let $Z\in \mathcal{P}$ and $Z\leq V$. Then by lemmas 3 and 11
we have $Z=Ann_V(Ann_V(Z))$. Hence $Z\in \mathcal{P}_V$.

b) It is trivial if $V=Z$.

Suppose $V\neq Z$; then $\{vzv : v\in Ann(V), z\in Z\}\neq \{0\}$.
Indeed, otherwise $\{vzv : v\in Ann(V), z\in Z\}=\{0\}$ and
$\{zv+vz : v\in Ann(V), z\in Z_+\}=\{0\}$. Hence, since for any
$z\in Z$ there exist $z_-$, $z_+\in Z_+$ such that $z=z_-+z_+$ we
have $\{zv+vz : v\in Ann(V), z\in Z\}=\{0\}$, i.e.  $Ann(V)\cdot
Z=\{0\}$ and $Z\subseteq Ann(Ann(V))$. This is impossible because
of $Ann(Ann(V))=V$. Let $B=\{vzv : v\in Ann(V), z\in Z\}$. We
assert that $Ann_Z(B)=V$. Suppose $Ann_Z(B)\neq V$; then there
exist $a\in Ann_Z(B)_+$ such that $a\notin V$ and $a\cdot
Ann(V)\neq \{0\}$. Hence there exists $v\in Ann(V)_+ $ such that
$a\cdot v\neq 0$. Suppose $avav=0$; then $avava=avcc^*va=0$, where
$a=cc^*$, $c\in A$. Hence $avc=0$ and $avcc^*=ava=0$. This is
impossible. Therefore $avav\neq 0$. By the definition $a\in Z$,
$vav\in B$ and $avav+vava\neq 0$. Note that $(avav)^*=vava$. At
the same time, by definition of $a$ $avav+vava=0$. This is a
contradiction. Therefore $Ann_Z(B)=V$. We have $B\subseteq
Ann_Z(V)$. Hence $Ann_Z(Ann_Z(V))=V$ and $Ann_Z(V)=\{vzv : v\in
Ann(V), z\in Z\}$. This concludes the proof of b).

c) is obvious. $\triangleright$

\medskip

An annihilator $V$ is said to be {\it modular}, if $\mathcal{P}_V$
is a modular lattice. The following lemma is valid by lemmata 8
and 13.

\medskip

{\bf Lemma 14.} {\it Let $A$ be a C$^*$-algebra and $V$ be an
Abelian  annihilator of $A$. Then

(a)  $\mathcal{P}_V$ is a boolean algebra,

(b)  each Abelian annihilator is modular. }

\medskip

The results of the given section can be summarized as the
following theorem.

{\bf Theorem 15.} {\it Let $A$ be a C$^*$-algebra and
$\mathcal{P}$ be the set of all annihilators of subsets of $A_+$.
Then

(a) $\mathcal{P}$ is a lattice with the order $\subseteq$,

(b) the annihilator $\{0\}$ is zero $\bf 0$ and $A$ is unit $\bf
1$ of the lattice $\mathcal{P}$,

(c) $\mathcal{P}$ is an ortholattice with the orthocomplementation
defined as $(\,\cdot\,)^\perp:\mathcal{P}\to \mathcal{P}$,
$(V)^\perp=Ann(V)$, $V\in \mathcal{P}$,

(d) If the C$^*$-algebra $A$ satisfies condition (A) then
$\mathcal{P}$ is an orthomodular lattice,

(e) elements $V$, $W\in \mathcal{P}$ are orthogonal as elements of
the ortholattice $\mathcal{P}$ if $V\cdot W=\{0\}$,

(f) the center of the ortholattice $\mathcal{P}$ coincides with
the set $Z(\mathcal{P})$ of all central annihilators of
$\mathcal{P}$,

(g) The lattice $\mathcal{P}$ is order complete,

(h)  the  center $Z(\mathcal{P})$ of $\mathcal{P}$ is a complete
boolean algebra.}

\medskip

{\bf Question 16.} {\it There arises the following question: When
does a C$^*$-algebra satisfy condition (A)?}

\medskip

{\bf Remark.} The lattice $\mathcal{P}$ of annihilators of a von
Neumann algebra $A$ is a sublattice of the lattice $J(A)$ of
$*$-weak closed inner ideals of the algebra $A$. The lattice
$J(A)$ is not orthomodular, but, since it possesses a
complementation, such concepts as orthogonality and center remain
meaningful nevertheless (see \cite{4}, \cite{5}). At the same
time, since $\mathcal{P}$ can be identified with the lattice
$P(A)$ of all projections in $A$, $\mathcal{P}$ is orthomodular.
Note that, in the case of anisotropic Jordan $*$-triples
annihilators in $\mathcal{P}$ are also inner ideals. In this case
elements of $\mathcal{P}$ are defined by Jordan multiplication.
Therefore the results in \cite{4} also hold for annihilators.

There exist many examples of uniformly closed two sided ideals of
a C$^*$-algebras which are not annihilators. Hence, since every
uniformly closed two sided ideal of a C$^*$-algebra is a
hereditary C$^*$-subalgebra no every hereditary C$^*$-subalgebra
is an annihilator.

\bigskip

\section{C$^*$-algebras of von Neumann type I}

\medskip

Now recall the definition of a C$^*$-algebra of type I. Let $A$ be
a C$^*$-algebra and $\pi:A\to B(H)$ be a representation of $A$,
where $H$ is a Hilbert space. The representation $\pi$ is said to
be {\it of type I}, if the von Neumann algebra, generated by $\pi
(A)$, is of type I. The C$^*$-algebra $A$ is said to be {\it of
type I}, if all representations of this algebra are of type I.

A C$^*$-algebra $A$ is called a CCR-{\it algebra}, if for every
representation $\pi:A\to B(H)$, where  $H$ is a Hilbert space,
such that $H\neq 0$ and the commutant of $\pi(A)$ in $B(H)$ is
${\bf C}1$ the operator $\pi(x)$ is a compact operator for every
$x\in A$.

A C$^*$-algebra $A$ is called a GCR-{\it algebra}, if each nonzero
factor-C$^*$-algebra of $A$ has a nonzero closed two sided
CCR-ideal. It is known that a C$^*$-algebra $A$ is a GCR-algebra
if and only if $A$ is a C$^*$-algebra of type I by Diximier
\cite{3}.

The theory of lattices of annihilators developed above allows us
to introduce the following definition.

{\it Definition.} A C$^*$-algebra $A$ is called a C$^*$-{\it
algebra of von Neumann type} I, if there exists an Abelian
annihilator $V$ in $\mathcal{P}$ such that $c(V)=A$.

\medskip

{\bf Proposition 17.} {\it Let $A$, $B$ be C$^*$-algebras, $\phi$
be a $*$-homomorphism of $A$ onto $B$. Then for every $S\subseteq
A_+$
$$
\phi(Ann(S))=Ann(\phi(S)).\,\,\,\,\,\,\,\,(**)
$$
}

{\it Proof.} If $A$, $B$ are von Neumann algebras and $\phi$ is
normal, then $A=ker \phi\oplus Ann(ker \phi)$ and $ker \phi=eA$,
$Ann(ker \phi)=(1-e)A$ for a central projection $e$ in $A$. Let
$S$ be a subset of $A_+$. Then $Ann(S)$ is a von Neumann algebra
and $Ann(S)=fAf$ for some projection $f\in A$ and
$$
\phi(Ann(S))=\phi((1-e)Ann(S))=\phi((1-e)fAf),
$$
$$
Ann(\phi(S))=Ann(\phi((1-e)S))=\phi(Ann((1-e)S))=\phi((1-e)fAf),
$$
since $\phi\vert_{(1-e)A}$ is a $*$-isomorphism of $(1-e)A$ onto
$B$. Hence the equality $(**)$ is valid.

Now, let $A$, $B$ be C$^*$-algebras, $\phi$ be a $*$-homomorphism
of $A$ onto $B$. Then by \cite[proposition 1.21.13]{Sak} $\phi$
has an extension to a normal $*$-isomorphism $\bar{\phi}$ of
$A^{**}$ onto $B^{**}$. Then $A^{**}=ker \bar{\phi}\oplus
Ann(\bar{\phi})$ and $ker \bar{\phi}=eA^{**}$, $Ann(ker
\bar{\phi})=(1-e)A^{**}$ for a central projection $e$ in $A^{**}$.
Note that, in this case $\bar{\phi}\vert_{(1-e)A}$ is a
$*$-isomorphism of $(1-e)A$ onto $B$.

Let $S$ be a subset of $A_+$. Then $a=ea+(1-e)a$,
$\phi(a)=\bar{\phi}(ea)+\bar{\phi}((1-e)a)=\bar{\phi}((1-e)a)$ for
every $a\in Ann(S)$. Therefore
$\phi(Ann(S))=\bar{\phi}((1-e)Ann(S))$. Similarly
$Ann(\phi(S))=Ann(\bar{\phi}((1-e)S))$. For every $b\in
Ann(\phi(S))$ there exists $c\in A$ such that
$b=\phi(c)=\bar{\phi}((1-e)c)$ and for every $v\in
\phi(S)=\bar{\phi}((1-e)S)$ we have $bv+vb=0$. Hence
$$
\bar{\phi}((1-e)c)\bar{\phi}((1-e)s)+\bar{\phi}((1-e)s)\bar{\phi}((1-e)c)=0
$$
and
$$
\bar{\phi}[((1-e)c)((1-e)s)+((1-e)s)((1-e)c)]=0
$$
for each $s\in S$ and
$$
Ann(\phi(S))\subseteq
\bar{\phi}(Ann_{(1-e)A}((1-e)S))=\bar{\phi}((1-e)Ann_{A}(S))=\phi(Ann(S)).
$$
Hence $Ann(\phi(S))\subseteq \phi(Ann(S))$.

Now, let $b$ be an element in $\phi(Ann(S))$. Then there exists
$c\in Ann(S)$ such that $b=\phi(c)$ and for every $s\in S$ we have
$cs+sc=0$. Hence
$$
\phi(cs+sc)=\phi(c)\phi(s)+\phi(s)\phi(c)=0
$$
and $bv+vb=0$ for every $v\in \phi(S)$. Therefore $b\in
Ann(\phi(S))$ and $\phi(Ann(S))\subseteq Ann(\phi(S))$. Thus
$\phi(Ann(S))=Ann(\phi(S))$. The proof is completed.
$\triangleright$

{\bf Theorem 18.} {\it  Let $A$ be a C$^*$-algebra on a Hilbert
space $H$, $w(A)$ be the weak closure of $A$ in $B(H)$ and $e$ be
an Abelian projection in $w(A)$ such that $A\cap ew(A)e\neq
\{0\}$. Then $Ann(Ann(A\cap ew(A)e))$ is an Abelian
C$^*$-algebra.}

{\it Proof.} Suppose $Ann(Ann(A\cap ew(A)e))$ is not Abelian. Let
$M=Ann(Ann(A\cap ew(A)e))$. Then for every factor representation
$\pi$ of $M$ we have
$$
\pi(M)=Ann_{\pi(M)}(Ann_{\pi(M)}(\pi(A\cap
ew(A)e))),\,\,\,\,\,\,\,(*)
$$
by proposition 17.

Since $A\cap ew(A)e$ is a hereditary Abelian C$^*$-subalgebra in
$M$ there exists a noncommutative factor representation $\pi$ of
type I of $M$, i.e. $w(\pi(M))$ is a noncommutative W$^*$-factor
of type I in $B(H_{\pi})$. We have $\pi(A\cap ew(A)e)\neq \{0\}$.
Indeed, otherwise $\pi(M)=\{0\}$ by $(*)$. It is clear that
$\pi(M)$ is not Abelian. Without loss of generality we may assume
$w(\pi(M))=B(H_{\pi})$.

Since $\pi(A\cap ew(A)e)$ is a hereditary Abelian C$^*$-subalgebra
in $\pi(M)$ there exists a projection $\bar{e}$ in $B(H_{\pi})$
such that
$$
w(\pi(A\cap ew(A)e))=\bar{e}B(H_{\pi})\bar{e}
$$
and $\bar{e}$ is Abelian. Suppose that $\bar{e}$ is a minimal
projection in $B(H_{\pi})$, i.e. $w(\pi(A\cap ew(A)e))$ is a
one-dimensional subspace in $B(H_{\pi})$ generated by $\bar{e}$.
Then
$$
w(\pi(A\cap ew(A)e))=\pi(A\cap ew(A)e)
$$
and $\bar{e}\in \pi(A\cap ew(A)e)$. Note that $\bar{e}\pi(A\cap
ew(A)e)\bar{e}=\pi(A\cap ew(A)e)={\Bbb C}\bar{e}$. Hence
$$
Ann_{\pi(M)}(Ann_{\pi(M)}(\pi(A\cap
ew(A)e)))=\bar{e}\pi(M)\bar{e}={\Bbb C}\bar{e}\neq \pi(M).
$$
The last inequality is a contradiction. Hence $w(\pi(A\cap
ew(A)e))$ is not one-dimensional, i.e. $\bar{e}B(H_{\pi})\bar{e}$
is not one-dimensional. Therefore $\bar{e}B(H_{\pi})\bar{e}$ is
not Abelian, i.e. $w(\pi(A\cap ew(A)e))$ is not Abelian. Therefore
$A\cap ew(A)e$ is not Abelian, but this is a contradiction. Thus
$M$ is Abelian. $\triangleright$

\medskip

{\bf Theorem 19.} {\it Let $A$ be a GCR-algebra on a Hilbert space
$H$. Then $A$ is a C$^*$-algebra of von Neumann type I. }

{\it Proof.} By lemma 4.4.4 in \cite{3} there exists a nonzero
element $x$ in $A$ such that $\pi(x)=0$ or $\pi(x)$ has rank 1 for
any representation $\pi$ of $A$. Hence
$\pi(xAx)=\pi(x)\pi(A)\pi(x)$ and $\pi(xAx)$ is commutative for
any representation $\pi$ of $A$. Therefore $xAx$ is a commutative
C$^*$-algebra.

Let $A_o$ be a maximal commutative $*$-subalgebra of $xAx$, then
for some maximal commutative $*$-subalgebra ${\bf A}_o$ of the
weak closure $w(xAx)$ of $xAx$ in $w(A)$ we have $A_o\subseteq
{\bf A}_o$. There exists a hyperstonian compact $Q$ such that
${\bf A}_o\cong C(Q)$. Let $e$ be the identity element of
$w(xAx)$. Then there exists a monotone increasing sequence $(x_n)$
(for example, an increasing approximate identity of $A_o$) such
that $\sup x_n=e$. Therefore the weak limit of the sequence
$(x_n)$ is $e$. Then $w(xAx)=ew(A)e$. Therefore by separately
weakly continuity of multiplication $ew(A)e$ is commutative.

We have $A\cap ew(A)e$ is an Abelian C$^*$-algebra. Let
$X=Ann(Ann(A\cap ew(A)e))$. Then by theorem 18 $X$ is an Abelian
C$^*$-algebra. Thus $A$ contains a nonzero Abelian annihilator
$X\in \mathcal{P}$.

Let $\{E_i\}$ be a maximal set of Abelian annihilators with
pairwise orthogonal central supports. We should prove that the
central support of $\bigvee_i E_i$ is $A$ that is if $c(E_i)$ is a
central support annihilator of $E_i$ for each $i$, then $\bigvee_i
c(E_i)=A$. If it is not true then $\bigvee_i c(E_i)<A$ and
$Ann(\bigvee_i c(E_i))\neq \{0\}$. Note that $Ann(\bigvee_i
c(E_i))$ is a central annihilator and a C$^*$-algebra. By theorem
4.3.5 in \cite{3} the annihilator $Ann(\bigvee_i c(E_i))$ is a
GCR-algebra. Hence there exists an Abelian annihilator $F$ in
$Ann(\bigvee_i c(E_i))$ with the central support $Z\subseteq
Ann(\bigvee_i c(E_i))$. This is contradicts the maximality of the
set $\{E_i\}$. Thus $\bigvee_i E_i=A$. Hence $A$ is a
C$^*$-algebra of von Neumann type I. $\triangleright$

{\bf Remark.} The converse of the statement of theorem 19 is not
true. For example, let $H_1$, $H_2$, $...$ be Hilbert spaces of
dimensions $1$, $2$, $\dots$ respectively. Then the C$^*$-algebra
$$
\sum_{n=1,2,\dots}^{\oplus}B(H_n)
$$
is not a GCR-algebra, but this algebra is a von Neumann algebra of
type I. Hence this algebra is a C$^*$-algebra of von Neumann type
I. Therefore the new class of C$^*$-algebras of von Neumann type I
is wider than the class of C$^*$-algebras of type I (that is the
class of GCR-algebras, \cite{3}).

\bigskip

\section{C$^*$-algebras without nonzero Abelian
annihilators}

\medskip

A C$^*$-algebra $A$ is called an NGCR-algebra, if this algebra
does not have nonzero two sided CCR-ideals.

\medskip

{\bf Theorem 20.} {\it Let $A$ be a C$^*$-algebra on a Hilbert
space $H$. Then $A$ is an NGCR-algebra if and only if $A$ does not
have a nonzero Abelian annihilator.}

{\it Proof.} Suppose $A$ is an NGCR-algebra and $A$ has a nonzero
Abelian annihilator $X$. By lemma 3 there exists a projection
$p\in A$ such that $w(X)=pw(A)p$. By separately weakly continuity
in $w(A)p$ of multiplication $w(X)$ is commutative. Hence $p$ is
Abelian. Let $e\in w(A)$ be the central support of $p$. Then
$ew(A)e$ is a von Neumann algebra of type I. We have $ew(A)e$ is a
von Neumann algebra of bounded linear operators on some Hilbert
subspace $H_o$ of $H$. Note that there exists a maximal
commutative subalgebra $A_o$ of $B(H_o)$ containing the
annihilator $X$. So $X\cap \mathcal{K}B(H_o)\neq \{0\}$, where
$\mathcal{K}B(H_o)$ is the algebra of all compact linear operators
on $H_o$. We have $I=ew(A)e\cap \mathcal{K}B(H_o)$ is a two sided
CCR-ideal of $ew(A)e$. Let $L=I\cap A$. Then $L\neq \{0\}$.
Indeed, $X\cap \mathcal{K}B(H_o)\subseteq I$ and $X\cap
\mathcal{K}B(H_o)\neq \{0\}$. By the definition $L$ is a two sided
CCR-ideal of $A$. This contradicts the assumption that $A$ is an
NGCR-algebra.

Suppose $A$ does not have a nonzero Abelian annihilator and there
exists a nonzero two sided CCR-ideal $I$ in $A$; then by the first
part of the proof of theorem 19 there exists an element $x\in I$
such that $xIx$ is a commutative C$^*$-algebra. We have
$xAx\subseteq I$. Therefore $xIx=xAx$ and by theorem 18 the
C$^*$-algebra $Ann(Ann(xAx))$ is commutative. This is a
contradiction of the supposition that $A$ does not have a nonzero
Abelian annihilator. Hence $A$ does not have a two sided
CCR-ideal. Therefore $A$ is an NGCR-algebra. $\triangleright$

Let $A$ be a C$^*$-algebra, $\mathcal{P}$ be the corresponding
lattice of annihilators. $\mathcal{P}$ is called {\it locally
modular}, if there exists a set $\{V_\xi\}$ of modular
annihilators with pairwise orthogonal central supports $\{Z_\xi\}$
such that $\sup_\xi Z_\xi=A$, in particular, if there exists  a
modular annihilator $V$ in $\mathcal{P}$ such that $c(V)=A$. It is
clear that in this case if $V=A$ then the lattice $\mathcal{P}$ is
modular. The lattice $\mathcal{P}$ is called {\it purely
nonmodular}, if there does not exist a nonzero modular annihilator
in $\mathcal{P}$. Recall that two annihilators $V$ and $W$ in
$\mathcal{P}$ are said to be {\it orthogonal}, if $V\cdot W=0$,
where $V\cdot W=\{vw:v\in V, w\in W\}$.

Let $\Xi$ be a set of indices and $\{Z_\xi\}_{\xi\in\Xi}$ be a set
of pairwise orthogonal central annihilators in $\mathcal{P}$. Let
$\sum^\oplus_{\xi\in\Xi} w(Z_\xi)$ be a set of subsets
$\{a_\xi\}_{\xi\in\Xi}$, where $a_\xi\in w(Z_\xi)$, with the
bounded set $\{\Vert a_\xi \Vert :\xi\in\Xi\}$.
$\sum^\oplus_{\xi\in\Xi} w(Z_\xi)$ is a von Neumann algebra with
the componentwise algebraic operations and the norm that is
defined as the least upper bound of the norms of the components
$a_\xi$.

Let $\mathcal{P}\vert_X=\{Y\in \mathcal{P}: Y\subseteq X\}$ and
$\mathcal{P}_X=\{Y\subseteq X: Ann_X(Ann_X(Y))=Y\}$, $X\in
\mathcal{P}$.

\medskip

{\bf Theorem 21.} {\it Let $A$ be a C$^*$-algebra on a Hilbert
space $H$, $w(B)$ be the weak closure of a subset $B\subseteq
B(H)$. Then there exist unique C$^*$-subalgebras $A_I$, $A_{II}$,
$A_{III}$ of $A$ such that

(a) $A_I$ is a C$^*$-algebra of von Neumann type I, there does not
exist a nonzero Abelian annihilator in the algebras $A_{II}$ and
$A_{III}$, the lattice ${\mathcal{P}_{A_{II}}}$ is locally
modular, the lattice ${\mathcal{P}_{A_{III}}}$ is purely
nonmodular.

(b) the C$^*$-subalgebras $A_I$, $A_{II}$, $A_{III}$ belong to
$Z(\mathcal{P})$,

(c) $A_I\oplus A_{II}\oplus A_{III}$ is a C$^*$-subalgebra of $A$
and
$$
Ann(A_I\oplus A_{II}\oplus A_{III})=\{0\}.
$$
}

{\it Proof.} Let $\{V_\xi\}_{\xi\in\Xi}$ be a maximal set of
Abelian annihilators with pairwise orthogonal central supports
$\{Z_\xi\}_{\xi\in\Xi}$, i.e. for any $\xi$ the annihilator
$Z_\xi$ is a central support of $V_\xi$ and $Z_\xi\cdot Z_\eta=0$
for every pair of different indices $\xi$ and $\eta$. Let
$\sum^\oplus_{\xi\in\Xi} V_\xi$ be a set of subsets $\{a_\xi:
\xi\in\Xi\}$, where $a_\xi\in V_\xi$, with the bounded set
$\{\Vert a_\xi \Vert: \xi\in\Xi\}$. The set
$\sum^\oplus_{\xi\in\Xi} V_\xi$ is a C$^*$-algebra with
componentwise algebraic operations and the norm, defined as
$$
\Vert a\Vert=\sup \{\Vert a_\xi \Vert:\xi\in\Xi\},
$$
where $a=\{a_\xi: \xi\in\Xi\}\in \sum^\oplus_{\xi\in\Xi} V_\xi$.

Indeed, the last assertion follows by $\sum^\oplus_{\xi\in\Xi}
V_\xi\subset \sum^\oplus_{\xi\in\Xi} w(V_\xi)$, where
$\sum^\oplus_{\xi\in\Xi} w(V_\xi)$ is a von Neumann algebra. By
separately weakly continuity of multiplication for every $\xi$
$w(V_\xi)$ is a commutative von Neumann algebra and there exists a
projection $p_\xi$ in $w(A)$ such that $w(V_\xi)=p_\xi
(w(A))p_\xi$. Let $p=\sum_{\xi\in\Xi} p_\xi$, where
$\sum_{\xi\in\Xi}$ is a weak limit of finite sums $\sum_{k=1}^m
p_k$, $\{p_k\}_{k=1}^m\subseteq \{p_\xi\}_{\xi\in \Xi}$. Then
$\sum^\oplus_{\xi\in\Xi} w(V_\xi)=\sum^\oplus_{\xi\in\Xi} p_\xi
(w(A))p_\xi= p(w(A))p$ and $p_\xi w(V_\xi)p_\xi$ is commutative
for all $\xi$. Hence $p(w(A))p$ is commutative and $p$ is an
Abelian projection. It is clear that $A\cap p(w(A))p$ is a
commutative C$^*$-subalgebra. Therefore by theorem 18 the
annihilator
$$
\sup_{\xi\in\Xi} V_\xi
$$
is a commutative C$^*$-algebra. Of course,
$A_I=\bigvee_{\xi\in\Xi} Z_\xi$ is the central support of
$\sup_{\xi\in\Xi} V_\xi$ and $A_I$ is a C$^*$-algebra of von
Neumann type I.

Let $Z=Ann(A_I)$. Then $Z$ is a central annihilator. Similarly we
can find a central annihilator $A_{II}$ in $\mathcal{P}_Z$ such
that $\mathcal{P}\vert_{A_{II}}$ is locally modular and
$w(A_{II})\oplus w(A_{III})\subseteq w(Z)$, where
$A_{III}=Ann_Z(A_{II})$, $Ann_Z(A_{II})$ is an annihilator of
$A_{II}$ in $Z$. By b) of lemma 13
$\mathcal{P}\vert_{A_{II}}=\mathcal{P}_{A_{II}}$ and by a) of
lemma 13 $A_{II}$, $A_{III}$ are central annihilators in
$\mathcal{P}$. By the definition of $A_{II}$ we have
$\mathcal{P}\vert_{A_{III}}$ is purely nonmodular. By b) of lemma
13 $\mathcal{P}\vert_{A_{III}}=\mathcal{P}_{A_{III}}$. It is clear
that $A_I\oplus A_{II}\oplus A_{III}$ is a C$^*$-subalgebra of
$A$. We have
$$
Ann(A_I\oplus A_{II}\oplus A_{III})=\{0\}.
$$
Uniqueness of $A_I$, $A_{II}$, $A_{III}$ is valid by their
definition. This concludes the proof. $\triangleright$

\medskip

{\bf Corollary 22.} {\it Let $A$ be a C$^*$-algebra on a Hilbert
space $H$, $w(B)$ be the weak closure of $B\subseteq B(H)$. Then
there exist unique C$^*$-subalgebras $A_I$, $A_{NGCR}$ of $A$ such
that

(a) the C$^*$-subalgebra $A_I$ is a C$^*$-algebra of von Neumann
type I, the C$^*$-subalgebra $A_{NGCR}$ is an NGCR-algebra,

(b) $A_I\oplus A_{NGCR}$ is a C$^*$-subalgebra of $A$ and
$$
Ann(A_I\oplus A_{NGCR})=\{0\}.
$$
}

{\it Proof.} The corollary follows by theorems 20 and 21.
$\triangleright$

\medskip

{\bf Remark.} By the theory developed above we can also introduce
notions of C$^*$-algebra of types II, III as follows: A
C$^*$-algebra $A$ is said to be of type II, if the lattice
${\mathcal{P}_A}$ is locally modular and there does not exist a
nonzero Abelian  annihilator in $A$. A C$^*$- algebra $A$ is said
to be of type III, if the lattice ${\mathcal{P}_A}$ is purely
nonmodular.

In the book of Dixmier "C$^*$-algebras and their representations"
\cite{3} the notion of a C$^*$-algebra of type I was introduced
and considered with the other equivalent notions as GCR-algebras
of Kaplansky and the notion of Makey. Then such notions as
representations of types II, III have also been introduced.
However the notions of C$^*$-algebras of type II and III have not
been introduced and investigated yet. The reason for it is that,
if a C$^*$-algebra has a representation of type II (of type III)
then this algebra necessarily has a representation of type III
(respectively of type II). Therefore it is impossible to introduce
the notions of C$^*$-algebras of types II and III using
representations of types II and III. An NGCR-algebra has
representations of types II and III, but does not have
representations of type I. As for the new notions, if a
C$^*$-algebra is of type II, then in this algebra does not exist a
nonzero central annihilator, being a C$^*$-algebra of type III or
I. Similarly, if a C$^*$-algebra is of type III, then in this
algebra does not exist a nonzero central annihilator, being a
C$^*$-algebra of type II or I.

Thus, theorem 21 is an analog of the type classification for
C$^*$-algebras.

\bigskip

\section{Classification of C$^*$-factors of von
Neumann type I}

\medskip

{\it Definition.} Recall that a C$^*$-algebra is called {\it a
C$^*$-factor}, if it does not have nonzero proper two-sided ideals
$I$ and $J$ such that $I·J=\{0\}$, where $I·J=\{ab :a\in I, b\in
J\}$. For example, every simple C$^*$-algebra is a C$^*$-factor.
Also, every W$^*$-factor is a C$^*$-factor.

{\bf Theorem 23.} {\it Let $A$ be a C$^*$-factor of von Neumann
type I on a Hilbert space $H$. Then $w(A)$ is a W$^*$-factor of
type I.}

{\it Proof.} Let $X$ be an Abelian annihilator in $\mathcal{P}$
such that $c(X)=A$. Then there exists an Abelian projection $e$
such that $w(X)=ew(A)e$ by lemma 3. Let $z$ be a central
projection in $w(A)$ such that $c(e)=z$, i.e. $z$ is a central
support of $e$. Then $X\subseteq zw(A)$. Let $I=zw(A)\cap A$.
Since $zw(A)$ is an ideal of $w(A)$, i.e. $zw(A)w(A)\subseteq
zw(A)$ we have $IA\subseteq I$. Hence $I$ is an ideal of $A$ and
$I=A$ since $A$ is simple. Hence $w(A)=zw(A)$ and $w(A)$ is of
type I.

Let $z$ be a central projection in $w(A)$ and $z<1$. Then $zX$ or
$(1-z)X$ is not equal to $\{0\}$. We note that $(1-z)e\neq 0$,
$ze\neq 0$ and $(1-z)e$, $ez\in w(X)$.

Let $Q$ be the topological space of multiplicative functionals on
$X$, $\bar{Q}$ be the topological space of multiplicative
functionals on $w(X)$. By lemma 9 $supp(Q)\subseteq supp(\bar{Q})$
and the set $supp(Q)$ is dense in $\bar{Q}$.

We have $V=\{t\in \bar{Q}: ze(t)\neq 0\}$ and $W=\{t\in \bar{Q}:
(1-z)e(t)\neq 0\}$ are close-open subsets of $\bar{Q}$ and
$supp(\bar{Q})=V\cup W$. Note that $V\cap supp(Q)$ is dense in $V$
and $W\cap supp(Q)$ is dense in $W$. Suppose $W$ does not contain
an open subset of $Q$; then $W\cap supp(Q)$ does not contain an
open subset of $Q$. In this case $V\cap supp(Q)$ is dense in $Q$.
Indeed, if $V\cap supp(Q)$ is not dense in $Q$ then the closure
$Cl(V\cap supp(Q))$ of $V\cap supp(Q)$ in $Q$ is not equal to $Q$,
i.e. $Cl(V\cap supp(Q))\neq supp(Q)$ and $supp(Q)\setminus
Cl(V\cap supp(Q))$ is an open subset in $Q$, that is contained in
$W\cap supp(Q)$. This is a contradiction.

Thus $V\cap supp(Q)$ is dense in $Q$. Hence every function $f$ in
the algebra $C(Q)$ of all real-valued continuous functions on the
locally compact space $Q$ is an unique extension of the function
$f_{V\cap supp(Q)}$ defined on $V\cap supp(Q)$. Therefore $C(Q)$
can be embedded in $C(V)=\{f\in C(\bar{Q}):\{x\in \bar{Q}:f(x)\neq
0\}\subseteq V\}$. Then every function $f$ in $C(V)$ has an unique
continuous extension on $\bar{Q}$ since $supp(Q)$ is dense in
$\bar{Q}$, that is $C(V)$ and $C(\bar{Q})$ can be identified in
the sense of $V\subseteq \bar{Q}$. Hence $V$ is dense in
$\bar{Q}$. Otherwise $supp(\bar{Q})\setminus Cl(V)$ is open and
nonempty in $\bar{Q}$. In this case $Cl(V)=V$ because $V$ is
close-open in $\bar{Q}$. Hence $W=supp(\bar{Q})\setminus Cl(V)$,
$supp(\bar{Q})\setminus Cl(V)$ is also close-open in $\bar{Q}$.
Then $C(\bar{Q})=C(V)\oplus C(W)$ and $C(W)\neq \{0\}$, that
contradicts the identifiability of $C(V)$ and $C(\bar{Q})$. Thus
$V$ is dense in $\bar{Q}$. Then $V=\bar{Q}$ since $V$ is a
close-open set in $\bar{Q}$. Hence $ze$ is an identity element of
$w(X)$, i.e. $ze=e$. Then by the previous part of the proof $z=1$.

Now, suppose $W$ contains an open subset $U$ of $Q$. Then the set
$C(U)$ of all functions $f\in C(Q)$ such that $\{x\in Q:f(x)\neq
0\}$ belongs to $U$ is a subalgebra of $C(Q)$. The algebra $C(U)$
has a nonzero function because $U\neq \oslash$ and $U$ is a
close-open set of $Q$. Therefore $(1-z)ew(X)\cap X$ and, hence,
$(1-z)w(A)\cap A$ are not empty. Hence $I_o=(1-z)w(A)\cap A$ is a
nonzero two sided ideal of $A$.

We have $w(I_o)$ is a nonzero two sided ideal of $w(A)$ and
$w(I_o)=w(A)$. Otherwise $Ann_{w(A)}(w(I_o))\neq \{0\}$ and
$Ann_A(I_o)\neq \{0\}$. Indeed, if $Ann_A(I_o)=\{0\}$ then
$Ann_{w(A)}(I_o)= \{0\}$ and $Ann_{w(A)}(w(I_o))= \{0\}$ by
separately weakly continuity of Jordan multiplication $a\circ
b=1/2(ab+ba)$. So, $Ann_A(I_o)\neq \{0\}$. Then $w(I_o)\oplus
w(Ann_A(I_o))\subseteq w(A)$ and $Ann_A(I_o)$ is a two sided ideal
of $A$ and
$$
I_o\cdot Ann_A(I_o)=\{0\}.
$$
Therefore $A$ is not a C$^*$-factor, what is impossible. Therefore
$w(I_o)=w(A)$. Then $1-z=1$. Hence $z=0$.

But $z$ is chosen arbitrarily. So $w(A)$ is a W$^*$-factor of type
I. $\triangleright$

\medskip

{\bf Theorem 24.} {\it Let $A$ be a simple C$^*$-algebra on a
Hilbert space $H$. Then $A$ is a CCR-algebra if and only if $A$ is
of von Neumann type I. }

{\it Proof.} By theorem 19, if $A$ is a CCR-algebra, then $A$ is
of von Neumann type I.

Now, suppose $A$ is of von Neumann type I. Let $\pi$ be a
representation of $A$ in a Hilbert space $H_{\pi}$. By the proof
of proposition 17 $\pi$ has an extension to a normal
$*$-representation $\bar{\pi}$ of $w(A)$ (the weak closure of $A$
in $B(H)$) onto $w(\pi(A))$ (the weak closure of $\pi(A)$ in
$B(H_{\pi})$). By theorem 23 $w(A)$ is a W$^*$-factor of type I.
We assert that $\pi(X)\neq \{0\}$ for every Abelian annihilator
$X\in \mathcal{P}_A$ such that $c(X)=A$. Indeed, there exists a
projection $p$ such that $w(X)=pw(A)p$. If $\bar{\pi}(p)=0$ then
$\bar{\pi}(q)=0$ for every projection $q$ in $w(A)$ equivalent to
$p$. Hence
$$
0=\sup_{q\sim p} \bar{\pi}(q)=\bar{\pi}(\sup_{q\sim p}
q)=\bar{\pi}(1_{w(A)})
$$
and $\bar{\pi}(w(A))=0$, where $1_{w(A)}$ is the identity element
of $w(A)$. This is a contradiction. Thus $\pi(X)\neq \{0\}$,
$X={\bf C}p$ and $\bar{\pi}(p)\neq 0$. Hence $\pi(A)$ is also a
simple C$^*$-algebra of type I and $w(\pi(A))$ is a W$^*$-factor
of type I. Since $\pi$ is chosen arbitrarily we have $A$ is a
CCR-algebra. $\triangleright$

{\it Definition.} Let $A$ be a C$^*$-algebra, $\mathcal{P}$ a
lattice of annihilators of $A$, $n$ be a cardinal number and $\Xi$
be a set of indices such that $\vert\Xi\vert=n$. We say $A$ is
{\it a C$^*$-algebra of type I$_n$}, if there is a set
$\{P_i\}_{i\in \Xi}$ of pairwise orthogonal Abelian annihilators
with the central support $A\in \mathcal{P}$ and $\sup_i
\{P_i\}_{i\in \Xi}=A$.

\medskip

{\bf Theorem 25.} {\it Let $A$ be a C$^*$-factor of von Neumann
type I. Then there exists a cardinal number $n$ such that $A$ is a
C$^*$-algebra of type I$_n$. }

{\it Proof.} Let $\{P_i\}$ be a maximal set of orthogonal Abelian
annihilators with a set of indices $\Xi$. It is clear that
$\mathcal{P}$ has only central elements $\{0\}$ and $A$. By
theorem 23 $w(A)$ is a W$^*$-factor of type I. By lemma 3 there
exists a projection $p_i\in w(A)$ such that $w(P_i)=p_iw(A)p_i$
for each $i$. By separately weakly continuity of multiplication
$w(P_i)$ is commutative for each $i$. Hence for every $i$ the
projection $p_i$ is Abelian and $\{p_i\}$ is an orthogonal set of
minimal projections in $w(A)$ (hence in the algebra $A$).

We suppose $\sup_i \{P_i\}_{i\in \Xi}<A$ in $\mathcal{P}$; then
$Ann(\sup_i \{P_i\}_{i\in \Xi})\neq \{0\}$. By theorem 24 $A$ is a
GCR-algebra. Hence $Ann(\sup_i \{P_i\}_{i\in \Xi})$ is also a
GCR-algebra and by the proof of theorem 19 there exists a nonzero
Abelian annihilator $X$ in $Ann(\sup_i \{P_i\}_{i\in \Xi})$. The
last statement contradicts maximality of $\{P_i\}$. Therefore
$\sup_i \{P_i\}_{i\in \Xi}=A$. Hence $A$ is a C$^*$-algebra of
type I$_n$, where $n=\vert\Xi\vert$. $\triangleright$

\medskip

By the arguments above a C$^*$-factor of type I can also be
defined as follows: a C$^*$-factor is said to be of type I if
there exists a nonzero Abelian annihilator in this factor.

\medskip

{\it Example.} Let $H$ be an infinite dimensional complex Hilbert
space. It is known that the space $\mathcal{K}B(H)$ of all compact
linear operators on $H$ is a simple C$^*$-algebra. Moreover it is
a CCR-algebra. The algebra $\mathcal{K}B(H)$ has a maximal
orthogonal set of minimal projections. Each of these minimal
projections generates an Abelian annihilator, which is isomorphic
to ${\bf C}$. These annihilators form a maximal orthogonal set of
Abelian annihilators in $\mathcal{P}_{\mathcal{K}B(H)}$ with the
central support $\mathcal{K}B(H)$. Hence $\mathcal{K}B(H)$ is a
simple C$^*$-algebra of type I$_n$, where $n=\dim(H)$.

\bigskip

\section{Equivalence relation in C$^*$-algebras}

\medskip

Let $A$ be a C$^*$-algebra on a Hilbert space $H$. Let $V$, $W$ be
annihilators of $\mathcal{P}$. We will write $V\approx W$, if
there exists a Banach subspace $B\subseteq A$ such that
$V_+=\{bb^*: b\in B\}$ and $W_+=\{b^*b: b\in B\}$.

\medskip

{\bf Lemma 26.} {\it Let $A$ be a C$^*$-algebra on a Hilbert space
$H$, $p$, $q$ be projections in $A$. Then $p\sim q$ if and only if
$pAp\approx qAq$. }

{\it Proof.} It is obvious that $pAp$, $qAq\in \mathcal{P}$.
Suppose $p\sim q$; then there exists an element $x\in A$ such that
$xx^*=p$, $x^*x=q$. Also we have $px=x$, $x^*p=x^*$, $xq=x$,
$qx^*=x^*$, $pAp=xAx^*$, $qAq=x^*Ax$. Let $B=\{xbx: b\in A\}$.
Then $B$ is a Banach space, $\{bb^*: b\in B\}=xA_+x^*$ and
$\{b^*b: b\in B\}=x^*A_+x$. Indeed, $xAx^*$ is a C$^*$-subalgebra
and $\{xbxx^*bx^*:b\in A\}\subseteq xA_+x^*$. Let $a\in xA_+x^*$.
Then there exists  $y\in xAx^*$ such that $a=yy^*$. We have
$y=xzx^*$ for some element $z\in A$. Also
$a=xzx^*xz^*x^*=xzx^*pxz^*x^*=x(zx^*)xx^*(xz^*)x^*$. Hence
$xzx^*x\in B$ and $xA_+x^*\subseteq \{bb^*:b\in B\}$. Therefore
$pAp\approx qAq$.

Conversely, suppose $pAp\approx qAq$ and let $B\subseteq A$ be a
corresponding Banach space; then there exists $b\in B$ such that
$bb^*=p$. Then $b^*b$ is a projection and $b^*b\in qAq$. Hence
$p\sim \leq q$. Similarly, there exists  $d\in B$ such that
$d^*d=q$. Then $dd^*$ is a projection and $dd^*\in pAp$. Hence
$p\geq \sim q$ and $p\sim q$ in $w(A)$, i.e. $p=xx^*$, $q=x^*x$
for some $x\in w(A)$. Hence $pw(A)p$ and $qw(A)q$ are isomorphic.
We have $pw(A)p_+=\{vv^*: v\in pw(A)q\}$, $qw(A)q_+=\{v^*v: v\in
pw(A)q\}$ and $B\subseteq pw(A)q$. The isomorphism is defined as
follows
$$
\phi: a\to x^*ax, a\in pw(A)p.
$$
In particular, if $a\in pw(A)p_+$ and $a=vv^*$, $v\in pw(A)q$ then
$$
\phi(vv^*)=v^*v, v\in pw(A)q
$$
and, if $a\in pAp_+$ and $a=vv^*$, $v\in B$ then
$$
\phi(vv^*)=v^*v, v\in B.
$$
Note that $B\subseteq pAq$ and, more precisely, $B=pAq$. Hence
$\phi\vert_{pAp}$ is an $*$-isomorphism of $pAp$ onto $qAq$. Since
$\phi(p)=q$ we have $p\sim q$ in $A$. $\triangleright$

\medskip

{\bf Lemma 27.} {\it Let $A$ be a C$^*$-algebra on a Hilbert space
$H$, $V$, $W$, $U$ be annihilators in $\mathcal{P}$, $p$, $q$ be
the identity elements of $w(V)$, $w(W)$ respectively  and
$S=pAq\cap A$. Then

(a) $V\approx W$ in $A$ if and only if $\{ss^*:s\in S\}=V_+$ and
$\{s^*s:s\in S\}=W_+$,

(b) if $V\approx W$ and $W\approx U$ in $A$, then $V\approx U$ in
$A$. }

{\it Proof.} We will prove (a). Let $V\approx W$ in $A$. Then
there exists a Banach space $B$ such that $\{bb^*:b\in B\}=V_+$
and $\{b^*b:b\in B\}=W_+$. It is easy to see, that $B\subseteq
pAq\cap A$. Hence $\{bb^*:b\in B\}\subseteq \{ss^*:s\in S\}$,
$\{b^*b:b\in B\}\subseteq \{s^*s:s\in S\}$. Since $\{ss^*:s\in
S\}\subseteq pAp\cap A_+$, $\{s^*s:s\in S\}\subseteq qAq\cap A_+$,
then $\{ss^*:s\in S\}=V_+$ and $\{s^*s:s\in S\}=W_+$. The converse
is obvious.

(b) Let $q$, $g$ be the identity elements of $w(W)$, $w(U)$
respectively, and $S=pAq\cap A$, $P=qAg\cap A$. Then by a) of
lemma 27 $V_+=\{aa^*: a\in S\}$, $W_+=\{a^*a: a\in S\}$,
$W_+=\{aa^*: a\in P\}$, $U_+=\{a^*a: a\in P\}$. We have
$$
\{(ab)(ab)^*: a\in S, b\in P\}=\{awa^*:a\in S, w\in W_+\}=
$$
$$
\{ac^*ca^*:a\in S, c\in S\}\supseteq \{(aa^*)^2:a\in
S\}=\{d^2:d\in V_+\}=V_+.
$$
At the same time, since $S=pAq\cap A$, $P=qAg\cap A$ we have
$$
\{(ab)(ab)^*: a\in S, b\in P\}\subseteq V_+.
$$
Hence $\{aa^*:a\in Q\}=V_+$, where $Q=\{ab: a\in S, b\in P\}$.

Similarly we get $\{a^*a:a\in Q\}=U_+$. $\triangleright$

\medskip

Note that by lemma 27 the relation $V\approx W$, where $V$, $W\in
\mathcal{P}$ is an equivalent relation of elements of the
ortholattice $\mathcal{P}$.

\bigskip

\section{C$^*$-factors without nonzero Abelian
annihilators}

\medskip

A projection $p$ in a C$^*$-algebra $A$ is said to be {\it
infinite}, if it is equivalent to a proper subprojection $q$ of
itself; and it is said to be {\it finite} otherwise. A simple
C$^*$-algebra is said to be {\it finite}, if every projection in
this algebra is finite.

It is known that there exists a faithful dimension function on any
modular lattice. Therefore, if for a given lattice does not exist
a faithful dimension function, then this lattice is not modular
\cite{16}.

\medskip

{\it Definition.} Let $A$ be a C$^*$-algebra, $\mathcal{P}$ be the
corresponding lattice of annihilators. $A$ is called a C$^*$-{\it
algebra of type} II, if $\mathcal{P}$ is locally modular. $A$ is
called a C$^*$-{\it algebra of type} II$_1$, if $\mathcal{P}$ is
modular. $A$ is called a C$^*$-{\it algebra of type} III, if
$\mathcal{P}$ is purely nonmodular.

\medskip

{\bf Theorem 28.} {\it A simple C$^*$-algebra of type II$_1$ is
finite.}

{\it Proof.} Let $A$ be a C$^*$-algebra of type II$_1$. Then the
lattice $\mathcal{P}$ of all annihilators in $A$ is modular. It is
known that a modular lattice is a continuous geometry. Then by the
results of von Neumann in \cite{16} there exists a faithful
dimension function $D$ on the lattice $\mathcal{P}$. By lemma 27
the values of the dimension function $D$ on equivalent
annihilators coincide. By lemma 26 and the additivity of $D$ every
projection in $A$ is finite. Hence $A$ is finite. $\triangleright$

A simple C$^*$-algebra $A$ is said to be {\it purely infinite} if
every nonzero hereditary subalgebra of $A$ contains an infinite
projection.

\medskip

{\bf Theorem 29.} {\it A simple purely infinite C$^*$-algebra is
of type III. }

{\it Proof.} Let $A$ be a simple purely infinite C$^*$-algebra. We
note that, each annihilator $X\in \mathcal{P}$ is a hereditary
C$^*$-subalgebra. Hence, by lemmas 2, 26 and 27 for each $X\in
\mathcal{P}$ there does not exist a nonzero faithful dimension
function on $X$. Hence each annihilator in $\mathcal{P}$ is
nonmodular. Hence, $A$ is of type III. $\triangleright$

\medskip

{\it Definition.} Let $A$ be a C$^*$-algebra, $\mathcal{P}$ be the
corresponding lattice of annihilators. $A$ is called a C$^*$-{\it
algebra of type} II$_\infty$, if $\mathcal{P}$ is locally modular
and the annihilator $A$ is nonmodular.

\medskip

{\it Example.} {\it 1.} Let $A$ be a W$^*$-factor of type
II$_\infty$ and $\{p_i\}$ be a maximal family of mutually
orthogonal equivalent finite projections in $A$ with $\sup_i
p_i=1$. Suppose $\{p_i\}$ is countable. Let
$$
\sum_{ij}^o p_iAp_j=\{\{a_{ij}\}:\,\,for\,\,each\,\,pair\,\,of\,\,
indices\,\,i,j,\, a_{ij}\in p_i Ap_j,\,and\,\,
$$
$$
\forall \varepsilon>0 \exists n_o\in N \,\,such\,\,that\,\,\forall
n\geq m\geq n_o \Vert
\sum_{i=m}^n[\sum_{k=1,\dots,i-1}(a_{ki}+a_{ik})+a_{ii}]\Vert<\varepsilon\}.
$$
Then $\sum_{ij}^o p_iAp_j$ is a C$^*$-algebra in relative to
componentwise algebraic operations, the bilinear operation and the
norm \cite{1}. Since $p_iAp_i$ is a simple finite C$^*$-algebra
for all $i$, then by the proof of theorem 8 in \cite{1} the
C$^*$-algebra $\sum_{ij}^o p_iAp_j$ is simple. Let
$\mathcal{A}=\sum_{ij}^o p_iAp_j$, $\mathcal{P}$ be the lattice of
all annihilators in $\mathcal{A}$. Then the annihilator
$\mathcal{A}\in \mathcal{P}$ is not modular. Indeed, let
$\{e_i\}$, $\{f_i\}$ be subfamilies of $\{p_i\}$ such that
$\vert\{e_i\}\vert=\vert\{f_i\}\vert$, $\{p_i\}=\{e_i\}\cup
\{f_i\}$ and $\{e_i\}\cap \{f_i\}=\{\oslash\}$. Let $e=\sup_i
e_i$. Then by the proof of lemma 26 $\{aa^*: a\in
e\mathcal{A}\}=e\mathcal{A}e$ and $\{aa^*:a\in
\mathcal{A}e\}=\mathcal{A}$. Hence $e\mathcal{A}e\sim \mathcal{A}$
in $\mathcal{P}$. Similarly $e\mathcal{A}e\sim f\mathcal{A}f$ and
$f\mathcal{A}f\sim \mathcal{A}$. Therefore there does not exist a
nonzero faithful dimension function on $\mathcal{A}$. So
$\mathcal{A}$ is not modular. Hence by the definition
$\mathcal{A}$ is a simple C$^*$-algebra of type II$_{\infty}$.

Suppose there exists an infinite projection $g$ in $\mathcal{A}$;
then $g=\sum_{ij}^o \{p_igp_j\}:=\lim_{n\to\infty} \sum_{i,j=1}^n
\{p_igp_j\}$ and $\Vert p_igp_i\Vert=1$ if $p_igp_i\neq 0$ for all
$i$. Since the projection $g$ is infinite then $\vert
\{i:p_igp_i\neq 0\}\vert=\infty$. Hence $\Vert
\sum_{i=1,\dots,n-1}(p_i gp_n+p_n gp_i)+p_n gp_n\Vert$ does not
converge to $0$ at $n\to \infty$. Hence $\{p_i gp_j\}\notin
\mathcal{A}$, i.e. $f\notin \mathcal{A}$. This is a contradiction.
Therefore $\mathcal{A}$ is finite.

Thus $\mathcal{A}$ is a finite simple C$^*$-algebra and of type
II$_{\infty}$.

{\it 2.} Let $A$ be a W$^*$-factor of type II$_\infty$, $\{p_i\}$
be a maximal orthogonal set of equivalent finite projections in
$A$ and $\sup_i p_i=1$. Let $\{\{p_j^i\}_j\}_i$ be the set of
infinite subsets of $\{p_i\}$ such that for all different indices
$\xi$ and $\eta$ $\{p_j^\xi\}_j\cap \{p_j^\eta\}_j=\oslash$,
$\vert\{p_j^\xi\}_j\vert=\vert \{p_j^\eta\}_j\vert$  and
$\{p_i\}=\cup_i \{p_j^i\}_j$. Let $q_i=\sup_j p_j^i$ for all $i$.
Then $\sup_i q_i=1$ and $\{q_i\}$ is a an orthogonal set of
equivalent projections. Then by theorem 9 in \cite{1} the
C$^*$-algebra $\sum_{ij}^o q_iAq_j$ is a C$^*$-factor with a
nonzero finite and an infinite projection. In this case
$\sum_{ij}^o q_iAq_j$ is not a von Neumann algebra. But
$\sum_{ij}^o q_iAq_j$ is a C$^*$-algebra of type II$_\infty$. This
assertion can be proved as in example 1 above.

{\bf Remark.} Let $(Fin)$ be the class of simple finite
C$^*$-algebras with no nonzero Abelian annihilators, $(PI)$ be the
class of simple purely infinite C$^*$-algebras, $(II_1)$,
$(II_\infty)$ and $(III)$ be the classes of simple C$^*$-algebras
of types II$_1$, II$_\infty$ and III respectively. Then by
theorems 28, 29 and examples above the following relations are
valid
$$
(II_1)\subset (Fin), (II_1)\neq (Fin), (Fin)\cap (II_\infty)\neq
\{\oslash\},(PI)\subset (III).
$$

\medskip

{\bf Theorem 30.} {\it For every C$^*$-factor $A$ one of the
following conditions holds:

    (a) $A$ is of type I$_n$, where $n$ is a natural number;

    (b) $A$ is of type I$_n$, where $n$ is an infinite cardinal number;

    (c) $A$ is of type II$_1$;

    (d) $A$ is of type II$_\infty$;

    (e) $A$ is of type III.
}

{\it Proof.} The theorem follows by theorems 21, 25 and by the
definitions of C$^*$-algebras of types II$_1$, II$_\infty$.
$\triangleright$

\medskip

{\bf Remark.} Note that in the case of von Neumann algebras the
definitions of C$^*$-algebras of types I$_n$, where $n$ is a
cardinal number, II$_1$, II$_\infty$ and III are equivalent to
(almost coincide with) the definitions of von Neumann algebras of
types I$_n$, II$_1$, II$_\infty$ and III respectively. By the
theory, developed above, there exist simple C$^*$-algebras of
types I$_n$, II$_1$, II$_\infty$ and III. At the same time there
exist only simple von Neumann algebras of types I$_n$, with $n$ is
finite, II$_1$, and III in the case of von Neumann algebras.

Note that, in the paper of R{\o}rdam \cite{15} it is given an
example of a simple C$^*$-algebra with a nonzero finite and an
infinite projection.

The approach to the classification problem for C$^*$-algebras
described in the given article may be closely connected to the
Elliott classification conjecture. Indeed, on the one hand,
theorem 30 is a completion of the theory, developed on the base of
the Elliott classification conjecture and other methods (see, in
particular,  \cite{2}, \cite{6}, \cite{7}, \cite{8}, \cite{9},
\cite{10}, \cite{11}, \cite{12}, \cite{13}, \cite{14}, \cite{15},
\cite{17}). On the other hand, the further developing the theory
based on the notions introduced and studied in the given article
may allow to add new type invariants to the list of the invariants
of the Elliott classification conjecture and form new
classification conjecture based on the Elliott classification
conjecture.

\medskip

\end{document}